\documentclass[journal]{IEEEtran}
\usepackage[english]{babel}
\usepackage{graphicx}
\usepackage{framed}
\usepackage[normalem]{ulem}
\usepackage{amsmath}
\usepackage{eurosym}
\usepackage{optidef}
\usepackage{cite}
\usepackage[pdf]{pstricks}
\usepackage{color, colortbl}
\definecolor{Gray}{gray}{0.9}
\usepackage{array}
\makeatletter
\newcommand{\thickhline}{\noalign {\ifnum 0=`}\fi \hrule height 1.2pt
    \futurelet \reserved@a \@xhline}
\newcolumntype{"}{@{\hskip\tabcolsep\vrule width 1pt\hskip\tabcolsep}}
\makeatother
\usepackage{amsthm}

\usepackage{amssymb}
\usepackage{amsfonts}
\usepackage{enumerate}
\usepackage{makecell}
\usepackage{soul}
\usepackage[shortlabels]{enumitem}
\usepackage{xfrac}
\usepackage[utf8]{inputenc}
\usepackage{nomencl}
\usepackage{multirow}

\ifCLASSINFOpdf
\usepackage{empheq}
\usepackage{xcolor}
\makeatletter
\usepackage[linesnumbered,ruled,vlined]{algorithm2e}
\usepackage{balance}

\newtheorem{problem}{Problem}
\newtheorem{remark}{Remark}
\allowdisplaybreaks

\makenomenclature
\usepackage{etoolbox}
\renewcommand\nomgroup[1]{%
  \item[\bfseries
  \ifstrequal{#1}{A}{Indices and Sets}{%
  \ifstrequal{#1}{B}{Parameters}{%
  \ifstrequal{#1}{C}{Optimization variables}{}}}%
]}

\begin{document}

\title{Demand Management for Peak to Average Ratio Minimization via Intraday Block Pricing}

\author{Carolina~Cortez,
         Andreas~Kasis,
         Dimitrios~Papadaskalopoulos
        and Stelios~Timotheou
\thanks{This work was funded by the European Union’s Horizon 2020 research and innovation program under grant agreements No. 891101 (SmarTher Grid) and No. 739551 (KIOS CoE), the project EnergyXchange - Enterprises/0618/0164 and from the Republic of Cyprus through the Directorate General for European Programs, Coordination, and Development.}
\thanks{C. Cortez, A. Kasis and  S. Timotheou are with the KIOS Research and Innovation Center of Excellence and the Department of Electrical and Computer Engineering, University of Cyprus, CY-1678 Nicosia, Cyprus (e-mail: cortez.carolina@ucy.ac.cy, kasis.andreas@ucy.ac.cy,  timotheou.stelios@ucy.ac.cy). 
{Dimitrios Papadaskalopoulos is with the Department of Electrical and Computer Engineering, University of Patras, Patras, Greece, 26504 (email: dimpap@upatras.gr).}}
\thanks{A preliminary version of this work has appeared in \cite{cortez2021}. Compared to \cite{cortez2021}, this paper considers the problem of designing the IBP prices with the objective to minimize the PAR of the power grid and also to ensure revenue adequacy and bill protection. In addition a solution approach, associated with the resulting optimization problem, has been developed and further discussion and simulations have been incorporated.}}

\maketitle

\begin{abstract}
Price based demand response schemes may significantly improve power system efficiency. Additionally, it is desired that such schemes yield improved power operation, by reducing the peak consumption. This paper proposes the \textit{Intraday Block Pricing} (IBP) scheme, aiming to promote effective demand response among consumers by charging their electricity usage based on intraday time-slots. To design the prices associated with the proposed scheme, we formulate a bilevel optimization problem that aims to minimize the Peak-to-Average Ratio (PAR) and simultaneously benefit the consumers and the utility company. The bilevel problem is converted into a single-level Mathematical Program with Equilibrium Constraints (MPEC). The resulting MPEC is non-convex and includes nonlinear constraints. Hence, to obtain a solution, it is relaxed into a Mixed Integer Linear Program by dealing with all nonlinearities. To evaluate the conservativeness of the proposed approach, a lower bound to the cost of the original bilevel problem is obtained. 
The applicability of the proposed scheme is demonstrated with simulations on various case studies, which exhibit a significant reduction in PAR and economic gains for the utility company and consumers. Moreover,  simulation results show that the solutions of the original and relaxed problems are equivalent, demonstrating the effectiveness of the proposed solution approach.
{Further simulation results demonstrate significant advantages in the performance  of the IBP scheme when compared to existing schemes in the literature.}
\end{abstract}


\renewcommand{\figurename}{Fig.{}}

\nomenclature[A, 01]{\(t\in\mathcal{T}\)}{Index and set of time-slots, $\mathcal{T}=\{1,\cdots,T\}$}
\nomenclature[A, 02]{\(c\in\mathcal{C}\)}{Index and set of clusters, $\mathcal{C}=\{1,\cdots,C\}$}
\nomenclature[A, 03]{\(f\in\mathcal{F}\)}{Index and set of number of consumption-blocks, $\mathcal{F}=\{1,\cdots,F\}$}
\nomenclature[B, 01]{\(n_c\)}{Number of consumers in cluster $c$}
\nomenclature[B, 02]{\(r\)}{Rate of return for the utility company ($r\geq1$)}
\nomenclature[B, 03]{\(\lambda^w_t\)}{Energy cost rate at time-slot $t$} 
\nomenclature[B, 04]{\(\check{q}, \hat{q}\)}{Lower and upper bounds for energy consumption breakpoints} 
\nomenclature[B, 05]{\(D_{tc}\)}{Baseline demand at time-slot $t$ for cluster $c$} 
\nomenclature[B, 06]{\(\tau_c\)}{Load shifting cost coefficient for cluster $c$} 
\nomenclature[B, 07]{\(\sigma_c\)}{Percentage of load shifting flexibility for cluster $c$}
\nomenclature[C, 01]{\(d^s_{tcf}\)}{Estimated demand response under the IBP scheme regarding time-slot $t$, cluster $c$, and consumption-block $f$} 
\nomenclature[C, 02]{\(d^o_{tcf}\)}{{Baseline demand at time-slot $t$ for cluster $c$ and  consumption-block $f$}} 
\nomenclature[C, 03]{\(d^{sh}_{tc}\)}{Demand shifted under the IBP scheme regarding time-slot $t$, and cluster $c$} \nomenclature[C, 04]{\(\xi\)}{Price increment for the IBP scheme}
\nomenclature[C, 05]{\(q_f\)}{Energy consumption breakpoint for the IBP scheme  consumption-block $f$} 
\nomenclature[C, 06]{\(\lambda^s_f\)}{Electricity price for the IBP scheme consumption-block $f$}
\nomenclature[C, 07]{\(d^{peak}\)}{Aggregated load peak of the resulting power market}

\printnomenclature

\section{Introduction}
\label{sec:introduction}

\textbf{Motivation:}
Demand response programs are expected to become widespread over the next years, particularly in the residential sector \cite{pinson2014}.
Although individual household loads are typically small, their large numbers result in a substantial aggregate consumption. In particular, 
the residential sector accounts for almost 30\% of the electricity market share of which nearly 80\% are currently being operated by regulated utility companies under concession agreements \cite{IEA}.
In such agreements, the regulator allows a company to provide power in a certain area. However, to avoid market power abuse, the regulator determines the  retail prices.
Hence, the design of effective market mechanisms that enable efficient demand response among households is highly important. Such schemes may also take into account additional regulatory rules that aim to improve the operation of the power network.

Price-based demand response programs use price signals to alter the behaviour of electricity consumers, enabling the time-shifting of electricity demand. Such schemes may be used to achieve a more balanced daily load shape, effectively reducing the peak-to-average ratio (PAR). The latter may yield multiple benefits, since a reduced PAR is associated with rarer load-shedding and  power outage  events, improved power quality, and more efficient use of the power system assets \cite{Strbac2008,pinson2014,soliman2014game,mohsenian2010autonomous}. The need for PAR reduction is expected to grow as more plug-in hybrid electric vehicles are incorporated in the power grid, since those are expected to significantly add to the  peak consumption \cite{mohsenian2010autonomous,gomez2003impact}. The appropriate utilization of demand-side resources through the curtailment of {the} PAR may also reduce the need for building expensive backup generators and help to accommodate further penetration of intermittent renewable generation.

Shifting flexible loads may also reduce household individual energy costs.
In particular, consumers may reduce their electricity bill by changing their consumption patterns in accordance with the power price signals.
Hence, in addition to improving the power system performance, price-based demand response programs can also benefit consumers \cite{pinson2014,ALLCOTT2011820}. 
Hence, the design of suitable price schemes that will simultaneously minimize the PAR and reduce the consumers' costs is highly important.

\textbf{Literature Review:}
Smart grid technologies are expected to  enable the widespread adoption of price-based programs through enhanced monitoring, communication, and control.
Among the proposed types of pricing schemes, those based on time-varying prices are largely discussed as potential solutions to promote demand response among consumers. For example, real-time pricing (RTP) can reflect the continuous variation in the cost of electricity supply, which is a desirable attribute that is attracting growing popularity \cite{namerikawa2015}.
{However, the application of RTP should factor in the diversity of the consumers' risk preferences. 
In particular, risk-averse consumers, 
which are mainly associated with residential consumers,
prefer to avoid the risk of facing a very high price at some periods and thus opt for a less dynamic pricing scheme, while 
risk-prone consumers  opt for RTP, opting to occasionally face very high prices in order to gain significant long-run savings.}
{Concerning the former type of consumers,} previous studies showed that continuous price changes might result in consumer fatigue, making them responsive to only extreme price changes \cite{deng2015,behrangrad2015}.

An alternative scheme, that yields more predictable price patterns, is called Time-of-Use (ToU). 
The ToU pricing scheme includes a set of different prices for different hours of the day, or even different seasons of the year, defined at the beginning of a tariff cycle (e.g., a year), that remains constant throughout the cycle. ToU schemes have 
demonstrated a shift in demand from periods of peak to periods of valley consumption \cite{Fahrioglu2000}.
A concern in such schemes is that they may result in the concentration of flexible loads at lower-priced periods, potentially causing larger load profile imbalances \cite{Papadaskalopoulos2016}.

Another price-based program is the Inclining Block Rate (IBR) scheme, where the electricity use is charged in consumption-blocks.
In such schemes, consumers that have higher electricity consumption throughout a charging period, typically one or two months, are subjected to higher prices.
This strategy has been used in different markets to address supply shortage issues \cite{Lin2013} or promote social equality \cite{hung2017}. However, IBR schemes do not necessarily reduce the PAR in power systems \cite{Zhou2018}. 
 {Moreover, a block of use scheme has been proposed in \cite{ma2019block} to enable fairer pricing among consumers.
 However, this work does not factor demand response in the design of prices 
 and does not aim to achieve PAR reduction.
 }

The use of price schemes as a means to reduce the PAR  has been considered in the literature. In particular, \cite{mohsenian2010autonomous,soliman2014game,fadlullah2013gtes,samadi2014real} propose demand-side management schemes that schedule the load profile of consumers through price signals to minimize the PAR.
However, the proposed schemes raise issues of scalability and transparency since consumers are required to engage with dynamic prices and communicate their scheduled load profiles.
Hence, an aim of this work is to design predictable prices for households that reduce individual electricity bills and minimize the PAR of the power grid. 

\textbf{Contribution:}
This paper proposes a novel pricing scheme, called the Intraday Block Pricing (IBP) scheme that aims to enable improved demand response. 
The IBP scheme mechanism resembles that of the IBR pricing scheme; however, instead of considering consumption-blocks on a monthly basis, it considers consumption-blocks in intraday time-slots (e.g. one or two hours).
The IBP scheme aims to yield higher prices at times of high consumption, motivating a shift of the demand to times of lower consumption. 
The selection of prices for the IBP scheme aims to minimize the PAR and simultaneously provide a benefit to the consumers participating, compared to a flat price scheme.

The IBP scheme encourages consumers to manage their demand within each time-slot, aiming to distribute the demand throughout the day. 
The price design enables the minimization of the PAR of the power grid and simultaneously ensures economic benefit for the consumers and revenue adequacy for the utility company. 
Since the prices directly depend on the demand at that time, the IBP scheme prevents the concentration of loads at some particular time period.
An important advantage of the IBP scheme is {that it provides a simple pricing structure, which can be easily understood by the consumers.
This is in contrast with more complex schemes, which require active price tracking to enable dynamic demand response,
that might be perceived as a highly complex engagement by residential consumers.
On the other hand, 
predictable price structures, such as those in the ToU, IBP and IBR schemes, reduce the required engagement from consumers 
{due to the simplicity of their perception.}
}

{Additionally, under the proposed IBP scheme, individual consumers are charged with higher electricity prices per consumption unit when their demand exceeds a given threshold and lower prices otherwise. This promotes a sense of fairness and social equity, since higher prices are paid by consumers with higher demand than those with lower consumption, such as low income families (see \cite{li20121} and also the case study considered in preliminary work \cite{cortez2021}).}

For the implementation of the IBP scheme, we propose an optimization-based approach that investigates, from the point of view of a regulated utility company, how the parameters of the IBP scheme should be efficiently designed in a residential market. 
To accommodate the consumers behaviour in our approach, an estimate of the demand response is integrated into the design of the IBP prices. {By adopting this method, consumers do not need to disclose any
real-time information and hence their habits cannot be accurately inferred, ensuring that their privacy is preserved.}

The hierarchical relationship between the design of the IBP prices (leader) and the demand response of residential consumers (followers) is modeled as a bilevel optimization problem. 
To solve this problem, it is converted into an equivalent single-level Mathematical Program  with  Equilibrium  Constraints (MPEC). The resulting MPEC is non-convex   and includes nonlinear constraints. Hence, to obtain a solution, it is relaxed into a Mixed Integer Linear Program (MILP).
Finally, a lower bound to the bilevel optimization problem is obtained to evaluate the conservativeness of the solution resulting from the MILP problem.

The practicality of the proposed approach is demonstrated with realistic simulations, which show that the presented scheme yields 
significantly reduced PAR and economic gains for the utility company and the consumers {and exhibits robustness to parametric uncertainty associated with the demand response model}.
In addition, the solution obtained from the relaxed MILP problem matches the estimated lower bound to the original non-convex bilevel problem.
{Moreover, simulation results demonstrate that the proposed IBP scheme yields improved robustness to parametric uncertainty compared to an optimized ToU scheme and also improved performance when the demand shifting flexibility is high.}
Hence, the simulation results demonstrate the applicability of the proposed scheme, {its advantageous properties compared to existing schemes in the literature and} the effectiveness of the solution approach.

In summary, the main contributions of this paper are:
\begin{itemize}
    \item The development of the IBP scheme  which aims to enable improved price-based demand response.
    \item The formulation of a suitable bilevel optimization problem for the efficient design of the IBP scheme parameters in a regulated residential market. 
    The problem aims to achieve PAR minimization and simultaneously ensure revenue adequacy for the utility company and economic benefit for the consumers.
    \item  The development of a solution approach for the considered non-convex bilevel problem, by relaxing it to a MILP problem. 
    To assess the conservativeness of the solution approach, an iterative procedure to obtain a lower bound to the original problem is also developed.
    Numerical results demonstrate that the relaxed MILP problem yields globally optimal solutions to the original problem.
    \item The provision of quantitative evidence for the beneficial impact of the IBP scheme on the residential consumers and the power system  {and for its advantageous properties compared to existing schemes in the literature,} using suitable case studies with realistic data.
\end{itemize}

\textbf{Paper structure:}
The remainder of this paper is structured as follows.  The proposed IBP scheme is presented in Section \ref{sec:IBP}. Section \ref{sec:model} describes and formulates the problem for the efficient selection of the IBP prices in a regulated residential market, while Section \ref{sec:solutionApproach} outlines the solution approach.
{Section \ref{sec_discussion} discusses possible extensions of the considered problem.}
Case studies that demonstrate the beneficial impact of the IBP scheme are discussed in Section \ref{sec:casestudies} and conclusions are drawn in Section \ref{sec:Conclusions}. {Finally, {Appendix A describes how network constraints may be incorporated in the problem formulation and Appendix B} includes the formulation of an optimal ToU scheme, used in the presented case studies.}

\section{Intraday Block Pricing Scheme}
\label{sec:IBP}
\begin{figure}
    \centering 
    \includegraphics[width=3in]{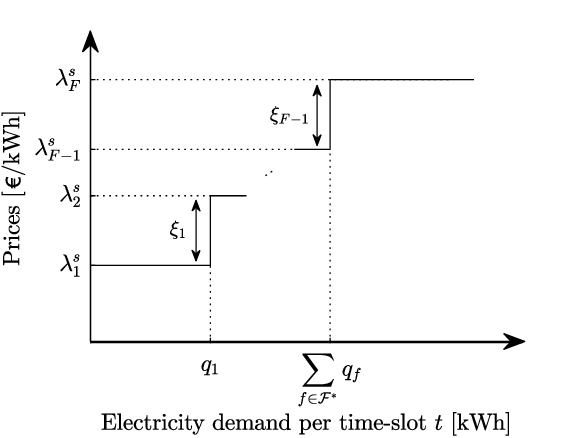}
    \vspace{-3mm}
    \caption{Schematic representation of the IBP scheme structure. The horizontal and vertical axes indicate the consumer electricity demand per pricing period and the IBP prices respectively.}
    \label{fig:IBPscheme} 
    \vspace{-5mm}
\end{figure}

The IBP scheme, illustrated in Fig. \ref{fig:IBPscheme}, consists of a \textit{price structure} comprised of the following parameters:
\begin{itemize}
    \item \textit{Time-slot -} the day is divided into $T$ time-slots of equal duration. 
    The set of time slots is denoted by $\mathcal{T}$.
    \item \textit{Clusters -} consumers with similar preferences are classified into clusters. Within each cluster $c\in\mathcal{C}$, where $\mathcal{C}$ is the set of clusters, there are $n_c$ consumers, such that $N\!=\!\sum_{c\in\mathcal{C}}n_c$ is the total number of consumers.
    \item \textit{Consumption-blocks -} the demand response to the IBP scheme for time-slot $t$ and cluster $c$, denoted by  $d^s_{tc}$, is split among consumption-blocks such that $d^s_{tc}=\sum_{f\in\mathcal{F}} d^s_{tcf}$, where $d^s_{tcf}$ is the demand response to the IBP scheme for time-slot $t$, cluster $c$, and consumption-block $f$, and $\mathcal{F}$  the set of consumption-blocks. We also define the set $\mathcal{F}^*=\{1,\cdots,F-1\}$, where\footnote{Note that for a set $\mathcal{A}$, we use $|\mathcal{A}|$ to denote its cardinality.}  $F = |\mathcal{F}|$.
    \item \textit{Electricity prices -}  consumption-blocks are associated with positive increasing price values satisfying $\lambda^s_{f+1}=\lambda^s_f+\xi_f$, $\forall(f\!\in\!\mathcal{F}^*)$, where  $\xi_f>0$ is the price increment from consumption-block $f$ to $f+1$, and $\lambda^s_f$ the price associated with consumption-block $f$. 
    \item \textit{Energy consumption breakpoints -} each consumption-block $f$ is limited by an amount of electricity demand $q_f>0$ such that $0\leq d^s_{tcf}\leq q_f$, $\forall(f\in\mathcal{F}^*)$ and $d^s_{tcF}\geq0$. The last consumption-block has no upper limit as the pricing scheme should not limit the total demand. 
    Note that since $\lambda_{f+1}>\lambda_{f}$, then  $d^s_{tc(f+1)}>0$ implies that $d^s_{tcf}=q_f$, i.e. for a non-zero consumption at block $f+1$, block $f$ must have been filled.
\end{itemize}

\section{IBP price structure design problem}
\label{sec:model}

In this section, we describe the considered problem, which aims to yield appropriate design parameters for the proposed IBP scheme. The latter is facilitated by formulating and solving  a suitable optimization problem.

\subsection{General Problem Description}
\label{sec:generalDescription}

We consider a smart power system with multiple residential consumers and a single utility company that, following the directives of a regulatory entity, designs the IBP prices. 
The utility company aims to replace an  existing flat pricing scheme with an (opt-in) IBP scheme {(\textit{Ass. 1)}\footnote{{For improved clarity, we enumerate the main assumptions behind the proposed approach throughout Section \ref{sec:generalDescription}.}}}.
{Similarly with flat pricing schemes, our solution approach designs the IBP price structure for a long period of time, 
aiming to include non-automated and traditionally less engaged consumers (e.g. ordinary households) in the demand response scheme.
Considering this time-frame, we opt to minimize the PAR 
aiming for a smoother power system operation.
The proposed IBP scheme could be employed simultaneously with dynamic schemes that may focus on alternative  power grid objectives  (e.g. power mismatch).
}

To increase the acceptance among consumers and to promote a smooth transition towards the IBP scheme, an incentive mechanism that  implicitly protects consumer bills is integrated into the design of the IBP prices. 
Revenue adequacy requirements are also considered to ensure the long-term sustainability of the utility company, i.e. the total revenue of the utility company must cover the cost of supply and allow a fair return. 
The desired minimum rate of return, denoted by $r$, is determined by the regulatory entity and it is typically based on the total value of the utility company's assets  (e.g., plant, equipment, working capital, and deductions for accumulated depreciation) {(\textit{Ass. 2})}.

{We consider exogenous parameters  that describe the wholesale market prices at different time slots, and assume that supply costs are equal to the wholesale market price times the procured energy {(\textit{Ass. 3})}.
The latter follows since the main residential sector supply costs are proportional to the energy supplied from the wholesale market, while alternative costs, such as  administrative
and operating costs, are significantly lower and usually fixed.
This approach is widely adopted in the literature  \cite{ferreira2013,emre20121,Mohsenian-Rad2010b}.
}

Typically, the energy cost rates (also known as wholesale prices) result from a strictly increasing, convex cost function {(\textit{Ass. 4})}. 
That is, the cost rates are high during peak demand due to the commitment of expensive generators, while during low demand periods more efficient generators are operating and thus the cost rates are lower. 

To ensure that consumers are motivated to act in a manner consistent with the rules established by the IBP scheme, the utility company integrates the consumers  demand response in the process of designing the IBP price structure. Towards this, we consider a communication channel from consumers to the utility company for billing purposes {(\textit{Ass. 5})}, which conforms with widely applied smart metering practices \cite{feng2019smart,akhavan2018power}. 
It should be noted that,  beyond billing purposes, household privacy is preserved, since no disclosure of information about the existence or assignment of local energy devices is required.
Based on this communication channel, the utility company can use the periodically collected data, associated with the flat prices, to estimate the load profile of consumers. In addition, consumers are assumed to react rationally to prices, aiming to maximize their welfare {(\textit{Ass. 6})} \cite{Bompard2007}. 

Although a large number of consumers can implement the proposed approach, 
{we assume that consumers can be grouped into a relatively small number of clusters that exhibit a well-defined behaviour  {(\textit{Ass. 7})}}.
The clustering method that could be adopted by the utility company is itself a widely discussed problem that is beyond the scope of this paper \cite{feng2019smart,Khafaf2021}. 
In addition, 
although individual demand is in general stochastic, 
we assume that consumer
clustering results in a predictable aggregate demand response {(\textit{Ass. 8})}.
{The latter is justified from the Central Limit Theorem (\cite[Ch. 2]{stroock2010probability1}) which suggests that the aggregation
of individual randomized consumption profiles tends to a normal distribution, when each cluster has a sufficient number of consumers.
Suitable clustering may yield a relatively small variance in the normal distribution describing each cluster, enabling an accurate description.}
{Furthermore, we opted for a single package IBP price structure to avoid discrimination among consumers, noting that a similar approach is adopted in several studies on retail pricing strategies \cite{momber20161,nguyen2016dynamic1}.}

To estimate the consumers response to the IBP scheme, we introduce an optimization problem that aims to minimize the electricity bill and load shifting inconvenience of consumers. 
The problem considers an aggregated, technology-agnostic demand response model {(\textit{Ass. 9})}, similar to \cite{ferreira2013,nojavan2017optimal,Sabita2013}, and assumes rational consumer behaviour. {In general, electricity demand exhibits flexibility potential on curtailment and shifting in time. 
However, the latter flexibility potential is much more prominent than the former, since, instead of simply avoiding using their loads at high price levels  consumers are more likely to shift the operation of their appliances to avoid higher prices, see e.g. \cite{Ye20181}. 
Hence, our approach focuses on encouraging a time shift in consumption and does not consider load curtailment {(\textit{Ass. 10})}.}

\begin{figure}
      \centering 
      \includegraphics[width=3.5in]{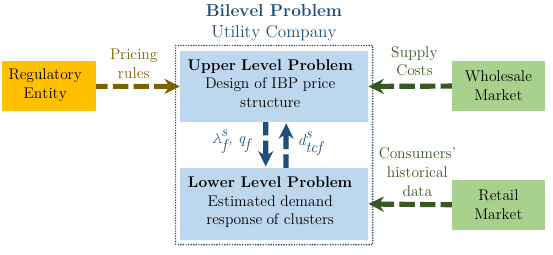}
      \vspace{-7mm}
      \caption{Bilevel optimization problem for the design of the IBP price structure.}
      \label{fig:bilevel}
      \vspace{-4mm}
\end{figure}

The utility company integrates the estimated demand response of consumers on its price-design problem. 
The hierarchical relationship between the design of the IBP price structure and the demand response of consumers  is modeled as a bilevel optimization problem, as illustrated in Fig. \ref{fig:bilevel}.

\subsection{Problem statement} 

In this section, we state the main problem considered in the paper associated with selecting suitable parameters for the IBP scheme by the utility company.

\begin{problem}\label{Problem_statement}
Select the parameters of the IBP price scheme to:
\begin{enumerate}[(R1)]
\item Minimize the PAR of the power grid.
\item Achieve revenue adequacy for the utility company. 
\item Protect consumer bills, ensuring no increase in individual bills due to the adoption of the IBP scheme. 
\item  Reflect the demand shifting flexibility of consumers.
\end{enumerate}
\end{problem}

Problem \ref{Problem_statement} aims to design the parameters of the IBP scheme such that Specifications (R1)-(R4) are satisfied.
In particular, it intends to minimize the PAR of the power system by promoting load shifting among households. In addition, it aims to achieve  revenue adequacy by ensuring that the monetary amount collected from the consumers covers retail costs and allows a fair return.
A further objective is to promote the  adoption of the IBP scheme by the consumers, through an incentive mechanism that implicitly protects their bills. 
Finally, an effective scheme should  integrate the estimated demand shifting flexibility of consumers in its price structure design.

\subsection{Utility company problem}
\label{subsec:UL}

The optimization problem considered by the utility company 
aims to satisfy Specifications (R1)-(R4), by suitable design of the IBP scheme parameters. The problem can be stated as

\begin{subequations}\allowdisplaybreaks\label{eq:UL}
    \begin{alignat}{2}
         \min_{\Xi_U}&                 \frac{\displaystyle\max_{t\in\mathcal{T}}\Bigg\{\sum_{c\in\mathcal{C}}\sum_{f\in\mathcal{F}}\!n_cd^s_{tcf}\Bigg\}}{\displaystyle\frac{1}{T}\sum_{t\in\mathcal{T}}\sum_{c\in\mathcal{C}}\sum_{f\in\mathcal{F}}n_cd^s_{tcf}} \span\span \label{eq:UL1}\\
         \text{s.t. }& \lambda^s_f=\lambda^s_1+(f-1)\xi, \span \forall(f\!\in\!\mathcal{F}), \label{eq:UL2}\\
        & \check{q}\leq q_f\leq\hat{q}, \span \forall(f\!\in\!\mathcal{F}^*), \label{eq:UL3}\\
        &  \sum_{t\in\mathcal{T}}\sum_{c\in\mathcal{C}}\sum_{f\in\mathcal{F}}\!\lambda^s_fn_cd^s_{tcf}\geq r\sum_{t\in\mathcal{T}}\sum_{c\in\mathcal{C}}\sum_{f\in\mathcal{F}}\lambda^w_tn_cd^s_{tcf}, \span\label{eq:UL4}\\
        &\sum_{t\in\mathcal{T}}\sum_{f\in\mathcal{F}}\lambda^s_fd^o_{tcf} \leq\lambda^o\sum_{t\in\mathcal{T}}D_{tc},
          &\forall(c\!\in\!\mathcal{C}),
		  \label{eq:UL5}
    \end{alignat}
\end{subequations}       
where $d^o_{tcf}$ represents the demand at time-slot $t$ for cluster $c$ associated with consumption-block $f\in\mathcal{F}$, given by

\begin{equation*}
    \begin{aligned}
        &d^o_{tcf}\! =\!\! 
        \min\!\!\left(\!\!\max\!\left(\!D_{tc}\! -\! \sum_{j=1}^{f-1}\! q_j,\! 0\!\right)\!\!, q_f\!\!\right)\!\!,
        \forall(t\!\in\!\mathcal{T},c\!\in\!\mathcal{C},f\!\in\!\mathcal{F^*}), \\
        &d^o_{tcF} = 
        \max\left(D_{tc} -\sum_{f\in\mathcal{F}^*}q_f, 0\right)\!\!,
        \forall(t\!\in\!\mathcal{T},c\!\in\!\mathcal{C}), 
    \end{aligned}
\end{equation*}  
which implies that $\sum_{f\in\mathcal{F}}d^o_{tcf}=D_{tc}, \forall(t\!\in\!\mathcal{T},c\!\in\!\mathcal{C})$.

In problem \eqref{eq:UL}, the set of decision variables is denoted by $\Xi_U\!=\!\{\xi;\ q_f,\ \forall f\in\mathcal{F}^*;\ \lambda^s_1\}$. 
In addition, $d^s_{tcf}$ denotes the estimated demand response of block $f$, cluster $c$ at time $t$, under the IBP scheme. These values are estimated by solving a lower level optimization problem, as explained in Section \ref{subsec:LL}.

Objective function \eqref{eq:UL1} aims to minimize the PAR, following Specification (R1). This is achieved by minimizing the peak (nominator of \eqref{eq:UL1}) over the average (denominator of \eqref{eq:UL1}) system demand over the considered time horizon $\mathcal{T}$.  

The price structure of the IBP scheme is captured through constraints \eqref{eq:UL2} and \eqref{eq:UL3}. To make the price structure more intuitive and increase the chances of its acceptance by the consumers, we consider fixed price increments $\xi$, such that $\xi_f = \xi$. We also consider that the energy consumption breakpoints have known lower and upper bounds $\check{q}$, $\hat{q}$, respectively. 
These bounds can be selected by the utility company or imposed by the regulatory entity.
It should be noted that the upper bound $\hat{q}$ should be selected to be sufficiently large to enable a broad set of solutions for \eqref{eq:UL}. A suitable option is to let $\hat{q}=\max_{t\in\mathcal{T},c\in\mathcal{C}}{D_{tc}}$.

To satisfy Specification (R2), constraint \eqref{eq:UL4} ensures that the total \textit{utility revenue} 
(left term) is \textit{at least as high as} the total \textit{utility cost} (right term) incurred within the time horizon $T$. We include a minimum rate of return $r\geq1$, determined by the regulator to incur a fair return to the utility company.   
As mentioned earlier, the supply cost is considered proportional to the energy supplied to the consumers, with cost rate $\lambda^w_t$ per time-slot $t$. {Our proposed approach considers exogenous data to reflect $\lambda^w_t$, under the assumption that the effect of the proposed scheme on wholesale prices is small. We opted to exclude the clearing process of the wholesale market to avoid shifting the focus of our study from the design of a suitable price structure.} 

An incentive mechanism for the adoption of the IBP scheme (Specification (R3)) is presented in  \eqref{eq:UL5}. The key idea is to ensure that the consumers electricity bill will not increase under the proposed scheme, compared to the existing flat price scheme, providing an incentive to opt-in the IBP scheme\footnote{Note also that problem \eqref{eq:LLs} ensures that the resulting demand from using the IBP scheme satisfies $\sum_{t\in\mathcal{T}}\sum_{f\in\mathcal{F}}\!\lambda^s_fd^s_{tcf} \leq \sum_{t\in\mathcal{T}}\sum_{f\in\mathcal{F}}\lambda^s_fd^o_{tcf}, \forall c \in \mathcal{C}$ since $d^s_{tcf}  = d^o_{tcf}, d^{sh}_{tc} = 0,  \forall(t\!\in\!\mathcal{T},c\!\in\!\mathcal{C},f\!\in\!\mathcal{F})$ is a feasible solution to it.
The latter guarantees that the implementation of the proposed scheme will not result in increased bills for consumers, regardless of whether they opt-in or not.}.
This is ensured by  \eqref{eq:UL5}, where $d^o_{tcf}$ describes the distribution of $D_{tc}$ at different consumption-blocks of the IBP scheme. 

\begin{remark}\label{Remark_existence}
The flat price $\lambda^o$ is selected to ensure that the market is budget balanced under a flat price and is given by 

        \vspace{-1mm}
\begin{equation} \label{eq:lambo}
        \lambda^o=r\sum\limits_{t\in\mathcal{T}}\sum\limits_{c\in\mathcal{C}}\lambda^w_t{n_c D_{tc}}\Big/\sum\limits_{t\in\mathcal{T}}\sum\limits_{c\in\mathcal{C}}n_c D_{tc}.
\end{equation}
Equation \eqref{eq:lambo} ensures that the revenue from selling energy at a flat price $\lambda^o$ equals the cost of the utility company. For consistency, parameter $r$ assumes the same value as in \eqref{eq:UL4}.
It should be noted that \eqref{eq:lambo} ensures the existence of a solution to \eqref{eq:UL}, by selecting the parameters of the IBP scheme to enable the original flat price.

\end{remark}

\subsection{Demand response model}
\label{subsec:LL}

The problem considered by the utility company to estimate the aggregated demand response of consumers, following Specification (R4), is presented below.

\begin{subequations}\allowdisplaybreaks\label{eq:LLs}
    \begin{alignat}{2}    
        \min_{\Xi_L}&
		   \sum_{c\in\mathcal{C}}\!n_c\left(\sum_{t\in\mathcal{T}}\sum_{f\in\mathcal{F}}\!\lambda^s_fd^s_{tcf}\!+\!\sum_{t\in\mathcal{T}}\!\frac{\tau_{c}}{2}{d^{sh}_{tc}}^2\right)\!   \span\label{eq:LLs1}\\
        \text{s.t. }& \sum_{f\in\mathcal{F}}d^s_{tcf}\!-\!d^{sh}_{tc}\!=\!D_{tc} 
          :\rho_{tc}, 
          &\forall(t\!\in\!\mathcal{T},c\!\in\!\mathcal{C}), \label{eq:LLs2}\\
        & 0\leq d^s_{tcf}
          :\mu^-_{tcf}, 
          &\forall(t\!\in\!\mathcal{T},c\!\in\!\mathcal{C},f\!\in\!\mathcal{F}), 
          \label{eq:LLs3}\\
        & d^s_{tcf}\!\leq q_f
          :\mu^+_{tcf}, 
          \quad\span\forall(t\!\in\!\mathcal{T},c\!\in\!\mathcal{C},f\!\in\!\mathcal{F}^*),
          \label{eq:LLs4}\\
        & -\sigma_cD_{tc}\leq d^{sh}_{tc}
          :\phi^-_{tc},  
          \span
          ~~\forall(t\!\in\!\mathcal{T},c\!\in\!\mathcal{C}), 
          \label{eq:LLs5}\\
        & d^{sh}_{tc}\leq \sigma_cD_{tc} 
          :\phi^+_{tc},  
          \span
          ~~\forall(t\!\in\!\mathcal{T},c\!\in\!\mathcal{C}), 
          \label{eq:LLs6}\\
        & \sum_{t\in\mathcal{T}}d^{sh}_{tc}=0 
          :\eta_c, 
          & \forall(c\!\in\!\mathcal{C}). \label{eq:LLs7}
    \end{alignat}
\end{subequations} 

In problem \eqref{eq:LLs}, the set of decision variables is denoted by $\Xi_L\!=\!\{d^{sh}_{tc},\ \forall(t\!\in\!\mathcal{T},c\!\in\!\mathcal{C});\ d^s_{tcf},\ \forall (t\!\in\!\mathcal{T},c\!\in\!\mathcal{C},f\!\in\!\mathcal{F})\}$.
Variable $d^{sh}_{tc}$ denotes the energy shifted in cluster $c \in \mathcal{C}$ at time-slot $t \in \mathcal{T}$, while variable $d^s_{tcf}$ denotes the estimated demand response of block $f$, cluster $c$ at time $t$, under the IBP scheme.

The objective function in \eqref{eq:LLs1} represents the collective cost perceived by consumers which includes the total \textit{consumer bill} paid to the utility company (first term) and the \textit{consumer shifting cost} due to load shifting in time (second term).
Parameter $\tau_{c}$ indicates the cost coefficient of load shifting at cluster $c$, describing the level of  discomfort that a consumer experiences when shifting his/her demand. 
{The cost associated with load shifting is described by a quadratic cost function that grows with the demand shifted (see e.g. \cite{62667201,erkoc2015game1}). 
It should be highlighted that the presented methodology and results can  be easily generalized to any convex cost function. 
} The objective function \eqref{eq:LLs1} inflicts a trade-off between the economic benefit provided by the IBP scheme prices $\lambda^s_f$ and the discomfort associated with the shifting cost coefficient $\tau_c$.

Constraints \eqref{eq:LLs2}-\eqref{eq:LLs7} express the time-shifting flexibility of different clusters using the IBP scheme.
Equality \eqref{eq:LLs2} ensures power balance at every time-slot $t\!\in\!\mathcal{T}$ and cluster $c\!\in\!\mathcal{C}$ such that the demand supplied by the utility company to the consumers, $\sum_{f\in\mathcal{F}}d^s_{tcf}$, minus the load shifting $d^{sh}_{tc}$ equals the baseline demand $D_{tc}$. The energy consumption partitioning among the consumption-blocks of the IBP scheme is expressed through constraints \eqref{eq:LLs3} and \eqref{eq:LLs4}. Constraints \eqref{eq:LLs5} and \eqref{eq:LLs6} provide lower and upper bounds on load shifting as a percentage $\sigma_c$ of the baseline demand.  Finally, constraint \eqref{eq:LLs7} ensures that the demand shifting is energy neutral within the considered time horizon.
Parameters $\sigma_c$ in constraints \eqref{eq:LLs5} and \eqref{eq:LLs6} quantify the amount of flexible demand that can be shifted  at each time-slot $t$.

Consumers can achieve economic benefit by shifting load from consumption-blocks with higher prices to consumption-blocks with lower prices.
Thus, the economic benefit results from the price increment $\xi$ between consumption-blocks, as in \eqref{eq:UL2}. 
Finally, the value of $q_f$ limits the amount of demand in each consumption-block $f$, thus it provides a guideline for the scheduling of load throughout the day.

It should be noted that \eqref{eq:LLs} is a strictly convex   problem. In addition,  it is intrinsically decoupled among consumption clusters, since consumer responses are independent between clusters.
Note also that the dual variables of problem \eqref{eq:LLs} are indicated after the colon of each constraint. For convenience we define the set of all dual variables $\Xi_D=\{\mu^{-}_{tcf},\ \forall(t\!\in\!\mathcal{T},c\!\in\!\mathcal{C},f\!\in\!\mathcal{F});\ \mu^{+}_{tcf},\ \forall(t\!\in\!\mathcal{T},c\!\in\!\mathcal{C},f\!\in\!\mathcal{F}^*);\ \phi^-_{tc},\ \phi^+_{tc},\ \rho_{tc},\ \forall(t\!\in\!\mathcal{T},c\!\in\!\mathcal{C});\ \eta_c,\ \forall(c\!\in\!\mathcal{C})\}$.

{The implementation of the proposed approach relies on the adopted model for flexible demand \eqref{eq:LLs}, and hence its parameters are of high importance.
Obtaining these parameters requires preliminary investigations with aim to estimate the consumer model parameters and associated behaviour.
These parameters may be obtained by pilot programs on demand response, consumer surveys and by making use of existing associated studies, see e.g. \cite{72327891, csereklyei2020price1,dergiades2008estimating1}. 
The procedure for obtaining these parameters, and generally the problem of accurate demand modelling, is of great importance but goes beyond the purposes of this study. In addition, assuming known parameters of consumers' response is not specifically associated with the proposed IBP scheme, but rather constitutes
a horizontal challenge for any type of pricing scheme explored in the relevant literature (e.g. RTP \cite{70736501,namerikawa2015,kim2016online1}, ToU \cite{ferreira2013,Datchanamoorthy20111,Zhou2018,62667201} and IBR \cite{Zhou2018,li20121} schemes).}
{Note also that any costs associated with implementing the IBP scheme are attributed to such preliminary studies and the potential installation of smart meters.}

{
\begin{remark}\label{rem_distrib_gen}
The presented demand response scheme does not include distributed generation. This can be incorporated trivially by replacing demand variables with net consumer demand variables, representing the difference between demand and generation, under the assumption that their values are non-negative at all times. Considering negative net demand values, which corresponds to generation, introduces several challenges due to the presence of reverse flows and induces further information requirements associated with the controllability of generation, and should be considered as part of future work.
\end{remark}
}

\subsection{Bilevel Optimization Problem for the IBP scheme}
\label{subsec:BL}

The resulting bilevel optimization problem is comprised of the upper level problem \eqref{eq:UL} and the lower level problem \eqref{eq:LLs}. The \textit{upper level problem}, describes the decision-making process of the utility company. For given estimated demand response of consumers, the  parameters of the IBP price structure are selected to satisfy Specifications (R1)-(R4). The \textit{lower level problem} estimates the demand response of the consumers to a given set of parameters for the IBP scheme.

The problem defined in  \eqref{eq:UL} and \eqref{eq:LLs} is challenging for two reasons. First, the problem is bilevel and hence cannot be handled by standard mathematical programming solvers. Second, the objective is fractional and also there are two non-convex constraints with bilinear terms (\eqref{eq:UL4} and \eqref{eq:UL5}). In the next section, we develop a solution approach to deal with these issues.  

{\begin{remark} \label{rem:networkconstraints}

The presented bilevel problem does not consider the presence of network constrains. 
We opted not to include such constraints in the presented problem formulation to avoid introducing additional complexity in an already mathematically dense paper and to
keep the focus on the structure and design process of the proposed IBP scheme.
In addition, it can be shown under practical conditions that the resulting demand response from the IBP scheme ensures that network constraints are not violated.
The above argument requires that the demand response from the originally imposed flat pricing scheme does not violate any network constraints.
The result can then be deduced using the results in \cite{dorfler2013synchronization1}, under realistic conditions on the network topology.
The intuition behind this argument follows by noting that the IBP scheme reduces both individual and aggregated peak consumption. 
\end{remark}
}

\section{Solution approach}
\label{sec:solutionApproach}

To solve the bilevel optimization problem introduced in the previous section we develop a two phase solution approach. First, we replace the lower level problem \eqref{eq:LLs}, by its equivalent Karush-Kuhn-Tucker (KKT) optimality conditions transforming the bilevel problem into a single-level Mathematical Program with Equilibrium Constraints (MPEC). Second, we relax the MPEC into a Mixed Integer Linear Program (MILP) by dealing with all nonlinearities.  

\subsection{Single level MPEC}
\label{subsec:MPEC}

To convert the bilevel problem \eqref{eq:UL}, \eqref{eq:LLs} into a single level MPEC, the lower level problem \eqref{eq:LLs} is replaced by its equivalent KKT conditions, presented in \eqref{eq:KKTLLs}. This transformation is allowed since \eqref{eq:LLs} includes a convex cost function and linear constraints  \cite{Boyd2004}.

\begin{subequations} \allowdisplaybreaks \label{eq:KKTLLs}
	\begin {alignat}{3}
        & n_c\tau_{c}d^{sh}_{tc}-\!\phi^-_{tc}+\phi^+_{tc}+\eta_{c}-\rho_{tc}=0, &\ \ \forall (t\!\in\!\mathcal{T},c\!\in\!\mathcal{C}),  \label{eq:KKTLLs1}\\
		& n_c\lambda^s_f-\mu^{-}_{tcf}+\mu^{+}_{tcf}+\rho_{tc}=0, \span \ \ \forall(t\!\in\!\mathcal{T},c\!\in\!\mathcal{C},f\!\in\!\mathcal{F}^*\!), \label{eq:KKTLLs2}\\
		& n_c\lambda^s_f-\mu^{-}_{tcf}+\rho_{tc}=0, \span \forall(t\!\in\!\mathcal{T},c\!\in\!\mathcal{C},f\!=\!F), \label{eq:KKTLLs2.1}\\
        &\mu^-_{tcf}d^s_{tcf}=0, \span \forall (t\!\in\!\mathcal{T},c\!\in\!\mathcal{C},f\!\in\!\mathcal{F}), \label{eq:KKTLLs3}\\
        &\!\mu^+_{tcf}\left(q_f\!-\!d^s_{tcf}\right)=0,
          \span\forall (t\!\in\!\mathcal{T},c\!\in\!\mathcal{C},f\!\in\!\mathcal{F}^*), \label{eq:KKTLLs4}\\
        &\phi^-_{tc}\left(d^{sh}_{tc}+\sigma_cD_{tc}\right)=0, 
          &\forall(t\!\in\!\mathcal{T},c\!\in\!\mathcal{C}), \label{eq:KKTLLs5}\\
        &\phi^+_{tc}\left(-d^{sh}_{tc}+\sigma_cD_{tc}\right)=0,
          &\forall(t\!\in\!\mathcal{T},c\!\in\!\mathcal{C}), \label{eq:KKTLLs6}\\
        &\mu^-_{tcf}\geq 0, \span \forall (t\!\in\!\mathcal{T},c\!\in\!\mathcal{C},f\!\in\!\mathcal{F}), \label{eq:KKTLLs7}\\
        &\mu^+_{tcf}\geq0,\ 
          \span\forall (t\!\in\!\mathcal{T},c\!\in\!\mathcal{C},f\!\in\!\mathcal{F}^*), \label{eq:KKTLLs8}\\
        &\phi^-_{tc}\geq 0, 
          &\forall(t\!\in\!\mathcal{T},c\!\in\!\mathcal{C}), \label{eq:KKTLLs9}\\
        &\phi^+_{tc}\geq 0,
          &\forall(t\!\in\!\mathcal{T},c\!\in\!\mathcal{C}), \label{eq:KKTLLs10}\\
        & \text{Constraints }\eqref{eq:LLs2}\text{-}\eqref{eq:LLs7}. \label{eq:KKTLLs11}
	\end{alignat}
\end{subequations}

Equalities \eqref{eq:KKTLLs1} and \eqref{eq:KKTLLs2}-\eqref{eq:KKTLLs2.1} are the optimality KKT conditions resulting from variables $d^{sh}_{tc}$ and $d^s_{tcf}$, respectively. Constraints \eqref{eq:KKTLLs3}-\eqref{eq:KKTLLs6} are the complementary slackness conditions. Constraints \eqref{eq:KKTLLs7}-\eqref{eq:KKTLLs10} are the dual feasibility conditions, while \eqref{eq:KKTLLs11} are the primal feasibility conditions \cite{Boyd2004}. Integrating the KKT conditions \eqref{eq:KKTLLs} into the bilevel problem yields the following MPEC: 

\begin{itemize}[label={}]
    \item Objective function: \eqref{eq:UL1} \hspace{3.3cm} (MPEC)
    \item Constraints: \eqref{eq:UL2}-\eqref{eq:UL5}, \eqref{eq:LLs2}-\eqref{eq:LLs7}, \eqref{eq:KKTLLs1}-\eqref{eq:KKTLLs10}.
    \item Variables: $\Xi_U$, $\Xi_L$ and $\Xi_D$.
\end{itemize}

\subsection{MILP Reformulation}
\label{sec:MILP}

The MPEC defined in Section \ref{subsec:MPEC} has three types of nonlinearities. 
The first nonlinearity type is associated with the presence of decision variables in the denominator of the objective function \eqref{eq:UL1}, yielding a fractional objective.
This fractional formulation can be transformed into a linear expression. 
In particular the lower level problem \eqref{eq:LLs} does not change the total demand of the households, since summing  \eqref{eq:LLs2} over $t\in\mathcal{T}$ and considering \eqref{eq:LLs7} yields  

\begin{align} \label{eq:D}
D=\sum_{t\in\mathcal{T}}\sum_{c\in\mathcal{C}}n_cD_{tc}=\sum_{t\in\mathcal{T}}\sum_{c\in\mathcal{C}}\sum_{f\in\mathcal{F}}n_cd^s_{tcf}    
\end{align}
where $D$ is a constant, allowing the denominator of the objective function \eqref{eq:UL1} to be removed. Hence, objective \eqref{eq:UL1} can be equivalently expressed as 

\begin{subequations}
    \begin{align}
         &\min d^{peak} \label{eq:nof}\\
         &d^{peak}\geq \sum_{c\in\mathcal{C}}\sum_{f\in\mathcal{F}}n_cd^s_{tcf}, \qquad \forall(t\!\in\!\mathcal{T}),\label{eq:dpeak}
    \end{align}
\end{subequations}
where $d^{peak}$ is a new variable.

The second type of nonlinearity is found in constraints \eqref{eq:UL4} and \eqref{eq:UL5} due to the products of variables $\lambda^s_f d^s_{tcf}$ and $\lambda^s_f d^o_{tcf}$, respectively.
To simplify the nonlinear term in \eqref{eq:UL4}, we substitute  \eqref{eq:UL2} into \eqref{eq:UL4} and consider \eqref{eq:D} yielding 

\begin{align}
     \label{eq:mani}
         \sum_{t\in\mathcal{T}}&\sum_{c\in\mathcal{C}}\sum_{f\in\mathcal{F}}n_c\lambda^s_fd^s_{tcf}= \nonumber\\
         &=\lambda^s_1D +\xi\sum_{t\in\mathcal{T}}\sum_{c\in\mathcal{C}}\sum_{f\in\mathcal{F}}n_c(f-1)d^s_{tcf}.
    \end{align}

Replacing \eqref{eq:mani} into \eqref{eq:UL4} and rearranging the terms yields

\begin{alignat}{2}
        &\lambda^s_1D+\sum_{t\in\mathcal{T}}\sum_{c\in\mathcal{C}}\sum_{f\in\mathcal{F}}(\xi f-\xi -r\lambda^w_t)n_cd^s_{tcf}\geq 0.\label{eq:new1}
\end{alignat}

Following a similar approach for constraint \eqref{eq:UL5} yields

\begin{alignat}{2}
        &\left(\lambda^s_1\!-\!\lambda^o\right)\sum_{t\in\mathcal{T}}D_{tc}\!+\!\xi\!\sum_{t\in\mathcal{T}}\sum_{f\in\mathcal{F}}(f-1)d^o_{tcf}\!\leq\!0,\  &\forall( c\in\mathcal{C}). 
		  \label{eq:new2}
\end{alignat}

All bilinear terms that appear in both \eqref{eq:new1} and \eqref{eq:new2} are functions of the variable $\xi$. Hence, our approach to deal with these nonlinearities is to fix $\xi$ and solve the problem over different values of $\xi \in [\check{\xi}, \hat{\xi}]$, where $\check{\xi}$ and $\hat{\xi}$ are sufficiently broad lower and upper bounds of the price increment, respectively. 
Consequently, $\xi$ is removed from the set of decision variables originating from the upper level problem.

The last category of nonlinearities involves the complementary slackness conditions \eqref{eq:KKTLLs3}-\eqref{eq:KKTLLs6} which are replaced with the respective linear constraints \eqref{eq:bigM}. Note that $M_1$, $M_2$, $M_3$, $M_4$ should be positive constants that are sufficiently large \cite{qiu2020}.

\begin{subequations} \allowdisplaybreaks\label{eq:bigM}
    \begin{alignat}{2}
        & d^s_{tcf}\leq \omega^1_{tcf}M_1,\quad \mu^-_{tcf}\leq(1-\omega^1_{tcf})M_1,\ \span\nonumber\\
        & & \qquad \omega^1_{tcf}\in\{0,1\},\ \forall(t\!\in\!\mathcal{T},c\!\in\!\mathcal{C},f\!\in\!\mathcal{F}), \label{eq:bigM1}\\
        & q_f\!-\!d^s_{tcf}\leq \omega^2_{tcf}M_2,\quad  \mu^+_{tcf}\leq(1-\omega^2_{tcf})M_2,\ \span\nonumber\\ 
        & &\omega^2_{tcf}\in\{0,1\},\ \forall(t\!\in\!\mathcal{T},c\!\in\!\mathcal{C},f\!\in\!\mathcal{F}^*), \label{eq:bigM2}\\
        & d^{sh}_{tc}+\sigma_cD_{tc}\leq \omega^3_{tc}M_3,\quad  \phi^-_{tc}\leq(1-\omega^3_{tc})M_3, \span\nonumber\\ 
        & &\omega^3_{tc}\in\{0,1\},\ \forall(t\!\in\!\mathcal{T},c\!\in\!\mathcal{C}), \label{eq:bigM3}\\
        & -d^{sh}_{tc}+\sigma_cD_{tc}\leq \omega^4_{tc}M_4,\quad  \phi^+_{tc}\leq(1-\omega^4_{tc})M_4, \span\nonumber\\ 
        & &\omega^4_{tc}\in\{0,1\},\ \forall(t\!\in\!\mathcal{T},c\!\in\!\mathcal{C}). \label{eq:bigM4}
    \end{alignat}
\end{subequations}

For convenience we define the set of all binary variables $\Xi_B=\{\omega^1_{tcf},\ \forall(t\!\in\!\mathcal{T},c\!\in\!\mathcal{C},f\!\in\!\mathcal{F});\ \omega^2_{tcf},\ \forall(t\!\in\!\mathcal{T},c\!\in\!\mathcal{C},f\!\in\!\mathcal{F}^*);\ \omega^3_{tc},\ \omega^4_{tc},\ \forall(t\!\in\!\mathcal{T},c\!\in\!\mathcal{C})\}$.
Hence, the MPEC problem is reformulated as a MILP problem as follows:

\begin{itemize}[label={}]
    \item Objective function: \eqref{eq:nof} \hspace{3.3cm} (MILP)
    \item Constraints: \eqref{eq:UL2}, \eqref{eq:UL3}, 
    \eqref{eq:LLs2}-\eqref{eq:LLs7}, \eqref{eq:KKTLLs1}-\eqref{eq:KKTLLs2}, \eqref{eq:KKTLLs7}-\eqref{eq:KKTLLs10}, \eqref{eq:dpeak}, \eqref{eq:new1}, \eqref{eq:new2} and \eqref{eq:bigM}.
    \item Variables: $\Xi_U\backslash\{\xi\}$, $\Xi_L$, $\Xi_D$, $\Xi_B$ and $d^{peak}$.
\end{itemize}

In the above formulation, the complementary slackness conditions are treated using \eqref{eq:bigM}, while variable $\xi$ is fixed. Hence, the resulting problem is a Mixed Integer Linear Program (MILP) which can be solved using standard optimization tools. 
The proposed MILP problem is a relaxed version of the MPEC problem, thus there is no guarantee that we can obtain a value for $\xi$ that globally minimizes the PAR. 
To evaluate the conservativeness of our solution approach,  we aim to compare the obtained solution with a theoretical lower bound of the global minimum of the original problem. The  approach to obtain this lower bound is described in the following section.

\subsection{Theoretical lower bound for the PAR minimization problem}
\label{sec:lowerbound}

In this section, we describe our approach to obtain a theoretical lower bound for the PAR minimization problem, described by \eqref{eq:UL}, \eqref{eq:LLs}.
In particular, we consider the load shifting flexibility associated with the lower level problem \eqref{eq:LLs}, and study load profiles that minimize the PAR but possibly violate \eqref{eq:UL4}-\eqref{eq:UL5}.

Our approach to obtain the lower bound, is by setting $|\mathcal{F}| = 2$ and letting $\xi$ be sufficiently large, such that consumers have a strong price incentive to limit their demand at $q_1$.
Note that, given a fixed value of $\xi$, the price values do not affect the resulting demand response, although these are important to satisfy constraints \eqref{eq:UL4}-\eqref{eq:UL5}. 
In addition, the case $|\mathcal{F}| > 2$ does not enable a reduced PAR compared to $|\mathcal{F}| = 2$. 
The latter follows since  when $\xi$ is selected to be sufficiently large, then consumers aim to keep their demand to no more than $q_1$ per time slot, provided this is feasible, i.e. does not violate \eqref{eq:LLs2}-\eqref{eq:LLs7}.
A potential addition of price blocks, when  $|\mathcal{F}| > 2$, would not hence provide additional incentive for shifting the peak demand, since all shifting ability is exhausted due to the large value of $\xi$.
Considering a larger set of price blocks could potentially result in improved demand distribution, and is associated with the satisfaction of \eqref{eq:UL4}-\eqref{eq:UL5}.

To obtain the value of $q_1$, we then consider an iterative procedure where, starting with $q_1^{(0)}=\check{q}$, the value of $q_1^{(i)}$ is increased by a small positive constant  $\varepsilon$  at each iteration $i$, such that $q_1^{(i)}=q_1^{(i-1)}+\varepsilon$. The process stops when $q_1^{(i)} \geq \hat{q}$ at some iteration $i$.
This approach enables two prices,  with the second being prohibitively expensive, encouraging  consumers to avoid requiring more than  $q_1$ energy units at some given time.
By considering a large number of values for $q_1 \in [\check{q}, \hat{q}]$, we consider the variability of this scheme to its changes.

The lower bound for the PAR is given by  $\underline{PAR}= \min_{{i\in\mathcal{I}}} \max_{t\in\mathcal{T}} \underline{d}^{s, (i)}_{t}$, where   $\underline{d}^{s, (i)}_{t}= \sum_{c\in\mathcal{C}}\sum_{f\in\mathcal{F}}{d^{s,(i)}_{tcf}}$,
{$\mathcal{I}$ is the set of iterations and ${d^{s,(i)}_{tcf}}$ the demand associated with time $t$, cluster $c$ and block $f$ at iteration $i$.}
Note the the lower bound to the PAR is associated with an iteration $i$ that minimizes the peak demand when $q_1 = q_1^{(i)}$.

{
\section{Potential Extensions}\label{sec_discussion}

In this section we discuss two potential extensions in the considered problem formulation.
In particular, we consider the incorporation of network constraints as well as uncertainty in the  demand response parameters and wholesale prices and discuss how they can be treated through the proposed solution methodology, presented in Section \ref{sec:solutionApproach}.
Although such extensions are relevant, we have opted not to incorporate them in the original problem formulation  for simplicity and to keep the paper focus on the structure and properties of the proposed IBP scheme.

\subsection{Network constraints}\label{sec_network_constraints}

The satisfaction of network constraints is important to ensure that the existing network infrastructure is able to satisfy the resulting demand from implementing the proposed IBP scheme.
In particular, network constraints require that the current network infrastructure suffices to satisfy the demand at all times.
An approach that demonstrates how such constraints can be included in the design of the IBP scheme prices, associated with three-phased balanced distribution networks, is presented in Appendix A.
Such constraints may be included in the design of the IBP prices in a relatively straightforward fashion since the existing solution approach, described in Section \ref{sec:solutionApproach}, directly extends to this case.
In addition, it is intuitive to note that, when radial distribution networks are considered,
 a sufficient condition for the network constraints to be satisfied under the IBP scheme is that those are not violated under flat pricing.
The latter follows since the IBP scheme results in a non-increase of individual  peak consumption, and therefore also in a non-increase of bus peak consumption.

\subsection{Uncertainty in model parameters}\label{sec_uncertainty}

A further potential extension in the considered problem  involves incorporating uncertainty in a set of model parameters.
Such parameters could be associated with the price demand response model or the wholesale prices which appear in the considered problem through the revenue adequacy constraint \eqref{eq:UL4}.
An approach to incorporate uncertainty in the design of IBP prices could entail:
\begin{enumerate}[label = (\roman*)]
    \item defining these parameters as random variables with known characteristics, e.g. drawn from specific probability distributions with known mean and variance, and
    \item defining the problem formulation as the minimization of the expected PAR value.
\end{enumerate}
A way to solve this problem is through\footnote{{A recent review on current methods  for power system optimization at the presence of uncertainty with suitable applications can be found in \cite{roald2023power}.}}
a scenario-based approach, where a set of values is drawn from the random variable distributions and then problem \eqref{eq:UL}, \eqref{eq:LLs} is simultaneously solved for all generated scenarios following the approach described in Section \ref{sec:solutionApproach}.
Such an approach would approximate the solution of the considered expected PAR minimization problem when a sufficiently large number of scenarios is included.
}

\section{Case studies}
\label{sec:casestudies}

This section considers realistic case studies  to evaluate the potential impact of the IBP scheme in minimizing the PAR and  benefit the consumers and the utility company. {In addition, it evaluates the robustness of the IBP scheme to parametric uncertainty associated with the demand response model and its scalability in terms of number of clusters.
Finally, it compares the performance of the IBP scheme with an optimized state-of-the-art scheme (ToU scheme).
}

\subsection{Simulation Setup} \label{subsec:setup}
\begin{figure}[t!]
    \centering 
    \includegraphics[width=3.3in]{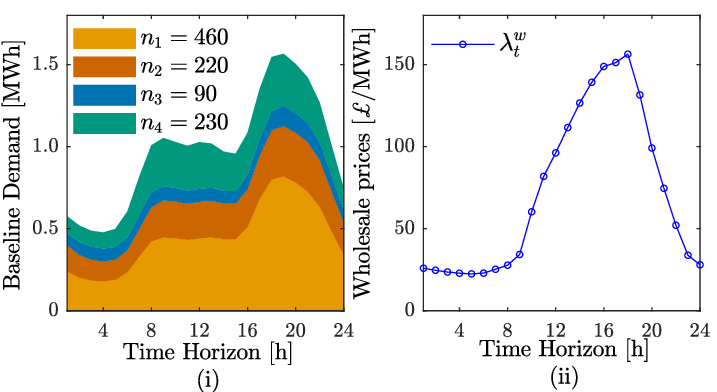}
        \vspace{-1mm}
    \caption{(i) Hourly baseline demand of each cluster. (ii) Wholesale prices.}
    \label{fig:input} 
    \vspace{-3mm}
\end{figure}

The examined case studies consider the IBP scheme in the context of a single day time horizon with hourly time-slots (i.e., $T\!=\!24$). The case study is built upon the power system of the United Kingdom (UK). The retail market consists of  $N\!=\!1,000$ residential consumers changing from flat pricing and opting into the proposed pricing scheme. The utility company categorizes its consumers into four clusters ($C\!=\!4$). The selected profiles and their respective contribution follow from \cite{Zoltan2016} which categorizes UK consumers into four archetypes named \textit{afternoon actives} ($c\!=\!1$), \textit{double risers} ($c\!=\!2$), \textit{winter spinners} ($c\!=\!3$) and  \textit{home lunchers} ($c\!=\!4$). The hourly baseline demand $D_{tc}$, $\forall(t\!\in\!\mathcal{T},c\!\in\!\mathcal{C})$ and the number of consumers per cluster $n_c$, $\forall(c\!\in\!\mathcal{C})$ are shown in Fig. \ref{fig:input}(i). The stacked area graphs indicate the aggregated load (summation of the individual demand of each consumer) with the contribution of each cluster represented by a different color. 

The assumed wholesale market prices follow the pattern of a typical winter day in the UK as illustrated in Fig. \ref{fig:input}(ii) \cite{qiu2020exploring}. The remaining operational costs associated with the utility company's activity are neglected; hence, it is assumed that $r\!=\!1$. 

The results from the optimal IBP price structure design approach proposed in Section \ref{sec:solutionApproach} are compared against a reference case with flat price $\lambda^o\!=\!0.083$ \pounds/kWh, calculated using  \eqref{eq:lambo} for the considered scenario. For simplicity, it is assumed that for each scenario the load shifting flexibility parameters and shifting costs are equal among clusters, that is $\sigma_c=\sigma$ and $\tau_c=\tau$ for all clusters $c\in\mathcal{C}$, respectively. 
The lower and upper bounds of the energy consumption breakpoints, such that $q_f\in[\check{q}, \hat{q}]$, were selected as follows
\vspace{-3mm}
\begin{equation}\label{eq:q}
    \check{q}=\min_{t\in\mathcal{T},c\in\mathcal{C}} D_{tc},
     \   \ 
    \hat{q}=\max_{t\in\mathcal{T},c\in\mathcal{C}} D_{tc}. \\
\end{equation}

Equation \eqref{eq:q} ensures that all consumption-blocks span through the baseline demand of different consumer clusters. Note that all optimization problems have been  solved using the mathematical programming software Gurobi \cite{Gurobi}.

\subsection{Impact of consumer shifting flexibility - $\sigma$}
\label{subsec:simul1}

This subsection examines the impact of the consumer shifting flexibility on the power grid PAR and the economic benefit perceived by the utility company and consumers. A total of 104 scenarios are considered for varying price increment values $\xi=\{0, 0.005, \cdots, 0.055, 0.06\}$ \pounds/kWh and different penetration of flexible loads $\sigma=\{10, 20, 30, 40\}$ \%, for two and three consumption-blocks ($|\mathcal{F}| \in\{2,3\}$). The implemented shifting cost is $\tau=0.03$ \pounds/kWh$^2$. Note that  the impact of $\tau$ is explored in Section \ref{subsec:simul2}.

\begin{figure}
    \centering
    \includegraphics[width=3.3in]{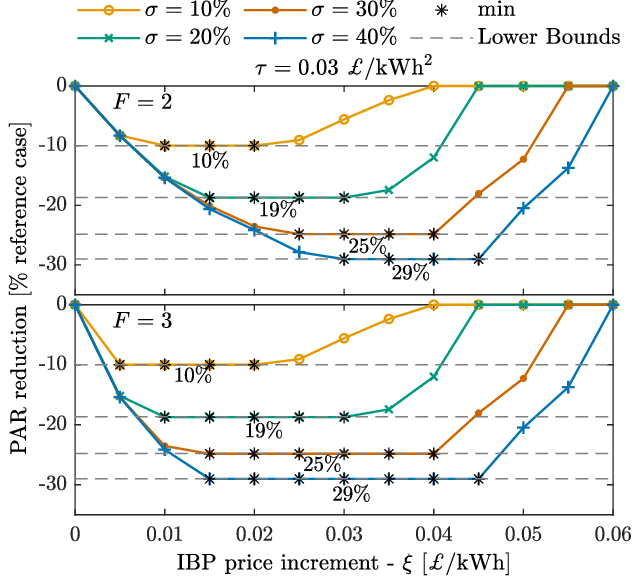}
    \vspace{-3mm}
    \caption{PAR reduction for increasing price increments $\xi$ and different load shifting flexibility values $\sigma$ when $F=2$ (top) and $F=3$ (bottom). \newline
    }
    \label{fig:PAR1}
    \vspace{-3mm}
\end{figure}

Figure \ref{fig:PAR1} depicts the PAR reduction as a percentage of the reference case for increasing values of the price increment $\xi$. Each plotted line stands for a different load shifting flexibility $\sigma$. 
Interestingly, all plots in Fig. \ref{fig:PAR1} have similar behaviour. Initially, the PAR does not present a significant reduction for small values of $\xi$, since the price increment  is not large enough to offer an economic incentive that outgrows the consumer shifting cost. 
Upon reaching a maximum reduction, the PAR remains constant and then decreases as $\xi$ increases. This occurs since the revenue adequacy and bill protection mechanisms, associated with constraints \eqref{eq:UL4} and \eqref{eq:UL5} respectively, limit the strategic choices of the utility company.

The theoretical lower bound for the PAR minimization problem, obtained following the approach presented in Section \ref{sec:lowerbound}, is depicted in Fig. \ref{fig:PAR1} with a dashed line. It can be observed that the solution from the relaxed MILP problem presented in Section \ref{sec:MILP} is equivalent to the optimal solution of the original bilevel problem. Note that, when the price increment $\xi$ becomes large, then the set of values that minimize the PAR is no longer feasible since the constraints \eqref{eq:UL4}-\eqref{eq:UL5} are violated.

The presented approach yields the set of values for $\xi$ for each given price structure and load shifting flexibility that enables the largest PAR reduction.
For example, if the market described in Section \ref{subsec:setup} has a load shifting flexibility of $\sigma=30\%$, then the utility company can potentially achieve a 25\% PAR reduction by selecting the parameters of the IBP price structure such that (i) $F=2$, $\lambda_1=0.079$ \pounds/kWh, $\xi=\{0.025, 0.03, 0.035, 0.04\}$ \pounds/kWh, and $q=1.119$ kWh or (ii) $F=3$, $\lambda_1=0.066$, $\xi=\{0.025, 0.03\}$ \pounds/kWh, and  $q=0.560$ kWh. 

\setlength{\tabcolsep}{4pt}
\begin{table}
    \caption{Impact of the price structure parameters on variations in economic and physical quantities for selected scenarios}
    \centering
     \begin{tabular}{l c c c c c c }
        \shortstack{\textbf{Scenarios}}
            & \multicolumn{2}{c}{\textbf{S1}} 
            & \multicolumn{2}{c}{\textbf{S2}}
            & \multicolumn{2}{c}{\textbf{S3}}\\[0.5pt]
        \thickhline
        \rowcolor{Gray} & & & & & &\\[-7.5pt]
        \rowcolor{Gray} \shortstack{$\sigma$ [\%]}
                            & \multicolumn{2}{c}{20} 
                            & \multicolumn{2}{c}{30} 
                            & \multicolumn{2}{c}{30}\\[0.5pt]
        & & & & & &\\[-7.5pt]
        \shortstack{$\tau$ [\pounds/kWh$^2$]}
            & \multicolumn{2}{c}{0.03} 
            & \multicolumn{2}{c}{0.03} 
            & \multicolumn{2}{c}{0.06}\\[0.5pt]
        \thickhline
        \rowcolor{Gray} & & & & & &\\[-7.5pt]
        \rowcolor{Gray} \shortstack{$\xi$ [\pounds/kWh]}
                            & \multicolumn{2}{c}{0.03} 
                            & \multicolumn{2}{c}{0.03} 
                            & \multicolumn{2}{c}{0.03}\\[0.5pt]
        & & & & & &\\[-7.5pt]
        \shortstack{$F$ [un.]} 
            & 2 & 3 & 2 & 3 & 2 & 3\\[1pt]
        \rowcolor{Gray} & & & & & &\\[-7.5pt]
        \rowcolor{Gray} \shortstack{$q_f$ [kWh]} 
                            & 1.145 & 0.501 & 1.119 & 0.689 & 1.118 & 0.533 \\[0.5pt]
        \rowcolor{Gray} & & 0.645 & & 0.430 & & 0.585 \\[0.5pt]
        & & & & & &\\[-7.5pt] 
        \shortstack{$\lambda^s_f$ [\pounds/kWh]}
            & 0.080 & 0.065 & 0.080 & 0.069 & 0.079 & 0.065 \\[0.5pt]
            & 0.110 & 0.095 & 0.110 & 0.099 & 0.109 & 0.095 \\[0.5pt]
            &  & 0.125 &  & 0.129 &  & 0.125 \\[0.5pt]
        \thickhline
        \rowcolor{Gray} & & & & & &\\[-7.5pt]
        \rowcolor{Gray} PAR [\%]$^{*}$ 
                            & 18.72 & 18.72 & 24.84 & 24.84 & 20.09 & 20.09\\[0.5pt]
        & & & & & &\\[-7.5pt]
        Utility cost [\%]$^{*}$ 
            & 4.14 & 4.14 & 5.88 & 5.96 & 5.07 & 5.07 \\[0.5pt]
        \rowcolor{Gray} & & & & & &\\[-7.5pt]
        \rowcolor{Gray} Consumer bill$^{\dagger}$ [\%]$^{*}$ 
            & 3.61 & 4.14 & 4.27 & 5.96 & 4.78 & 4.97 \\[0.5pt]
        & & & & & &\\[-7.5pt]
        Consumer total cost [\%]$^{*}$ 
            & 2.77 & 3.3 & 2.74 & 4.43 & 2.63 & 2.83 \\[0.5pt]
        \thickhline
        \multicolumn{7}{l}{} \\[-5pt]

        \multicolumn{7}{l}{\text{$^{*}$: Reduction in percentage compared to the reference case}}\\[1pt]
        \multicolumn{7}{l}{\text{$^{\dagger}$: Equivalent to utility revenue}}
    \end{tabular}
    \label{tab:simul1}
    \vspace{-3mm}
\end{table}

Table \ref{tab:simul1} examines the interaction between the price structure design and the demand response of consumers for three different scenarios (S1-S3) for two ($F=2$) and three ($F=3$) consumption-blocks. Table rows 1-2, 3-6 and 7-10 depict the consumer parameters, the optimal IBP price structure, and the performance metrics of each scenario respectively. 
Note that Scenario S2 is similar to S1 but has increased load shifting flexibility ($\sigma=30\%$ compared to $20\%$), while Scenario S3 is similar to S2 but has increased shifting cost  ($\tau=0.06$ compared to $\tau=0.03$). 
Comparing Scenarios S1 and S2 shows that a 50\% increase in load shifting flexibility results in approximately 40\% reduction in terms of utility cost and consumer bills. The same comparison for Scenarios S2 and S3 yields a 16\% performance reduction in both metrics. 
It should be noted that cases $F=2$ and $F=3$ yield the same maximum PAR reduction, which is in line with the discussion in Section \ref{sec:lowerbound}. 
However, the additional flexibility offered by having the extra pricing block may result in improved utility and consumer  costs, as demonstrated in Table \ref{tab:simul1}.

Considering the different metrics across all scenarios yields three important observations. First,  the PAR is reduced by more than 15\% compared to the reference case. 
Second, the utility operating cost, as reflected by the  consumer bills, is improved by 3-6\%. Finally, the results indicate that the reduction of the  consumer bills outgrows the perceived increase on shifting cost; hence, the consumer experience is improved with the implementation of the IBP scheme.  

\begin{figure}
    \centering
    \includegraphics[width=3.55in]{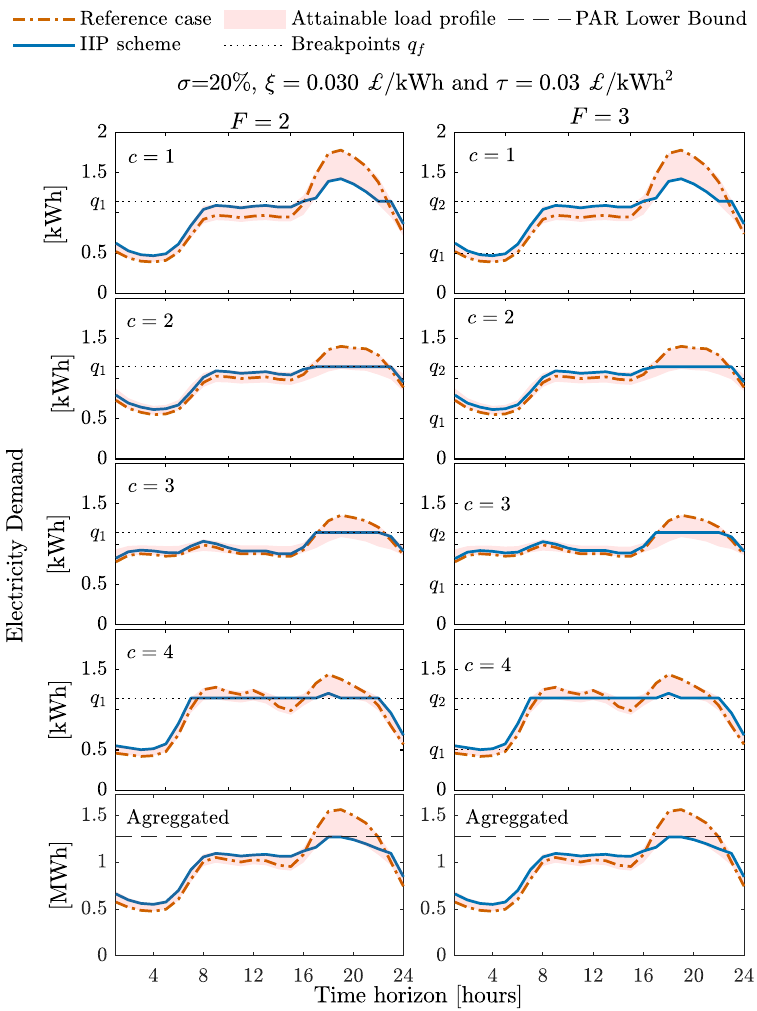}
    \vspace{-4mm}
    \caption{Electricity demand profile under Scenario S1. The red dashed and blue continuous lines depict the reference case, i.e. $D_{tc}$, and the demand response to the IBP scheme after load shifting, i.e. $\sum_{f\in\mathcal{F}}d^s_{tcf}$.}
    \label{fig:simul1}
    \vspace{-3mm}
\end{figure}

To illustrate how consumers respond to the IBP scheme, Fig. \ref{fig:simul1} provides the hourly demand of each cluster $c\in\mathcal{C}=\{1,2,3,4\}$, and the aggregated load profile experienced by the grid under Scenario S1 for $F=2$ and $F=3$. The load shifting flexibility is $\sigma=20\%$. The attainable load profile for each cluster $c\in\mathcal{C}$, indicated by the red shaded areas, are bounded from below and above by $\underline{d}^s_{tc}=\min_{i\in\mathcal{I}}\sum_{f\in\mathcal{F}}{d^{s,(i)}_{tcf}}$ and $\overline{d}^s_{tc}=\max_{i\in\mathcal{I}}\sum_{f\in\mathcal{F}}{d^{s,(i)}_{tcf}}$, respectively. The aggregated attainable load profiles are bounded by $\underline{d}^s_{t}=\min_{i\in\mathcal{I}}\sum_{c\in\mathcal{C}}\sum_{f\in\mathcal{F}}d^{s,(i)}_{tcf}$ and
  $\overline{d}^s_{t}=\max_{i\in\mathcal{I}}\sum_{c\in\mathcal{C}}\sum_{f\in\mathcal{F}}d^{s,(i)}_{tcf}$.
Note that the iterations associated with the derivations of these  bounds, and similarly the definition of variables $d^{s,(i)}_{tcf}$ follow from the approach presented in Section \ref{sec:lowerbound}.

From Fig. \ref{fig:simul1} it can be observed that the utility company strategically selects the breakpoint values, highlighted in the black dotted lines, to encourage consumers to monitor their consumption within each pricing period and shift their flexible load from periods of high to periods of low consumption. These strategic pricing signals reduce the PAR and lead to decreased utility operating costs and consumer bills. 
   
\subsection{Impact of shifting cost - $\tau$} 
\label{subsec:simul2}

This section analyzes the impact of the consumer shifting cost. 
A total of 156 scenarios are considered for varying price increments $\xi=\{0, 0.005, \cdots, 0.055, 0.06\}$ \pounds/kWh and different values of shifting costs $\tau=\{0,0.02,\cdots,0.08,0.1\}$ \pounds/kWh$^2$, for two and three consumption-blocks ($|\mathcal{F}|\in\{2,3\}$). The load shifting flexibility is set to $\sigma=30\%$.

\begin{figure}
    \centering
    \includegraphics[width=3.3in]{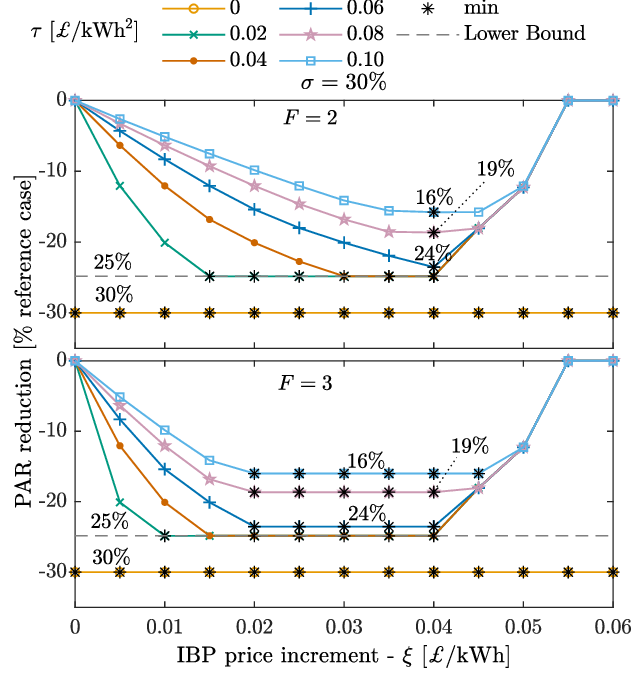}
    \vspace{-5mm}
    \caption{PAR reduction as a function of the price increment $\xi$ for different values of the load shifting cost $\tau$ when $F=2$ (top) and $F=3$ (bottom).}
    \label{fig:PAR2}
    \vspace{-3mm}
\end{figure}

Figure \ref{fig:PAR2} depicts the PAR reduction as a percentage of the reference case for increasing values of the price increment value $\xi$.
Each plotted line illustrates the PAR reduction for different shifting cost coefficients $\tau$. Similar to Fig. \ref{fig:PAR1}, the observed PAR reduction is a unimodal function of $\xi$. In this case, however, the PAR reduction decreases from 30\% to 16\% as the shifting cost coefficient increases from $\tau=0$ to $\tau=0.1$  \pounds/kWh$^2$. 
The latter follows since the willingness of consumers to shift their demand reduces as the  shifting cost increases, although  their shifting flexibility remains the same. 

\begin{figure}
    \centering
     \includegraphics[width=3.5in]{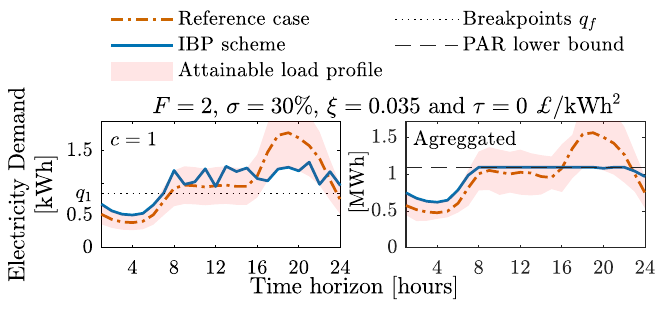}
     \vspace{-5mm}
    \caption{Electricity demand profile when consumers have zero shifting cost ($\tau=0$ \pounds/kWh$^2$) for cluster 1 (left) and all consumer clusters (right).}
    \label{fig:simul2}
    \vspace{-3mm}
\end{figure}

An interesting case depicted in Fig. \ref{fig:PAR2} is when $\tau=0$, i.e. when there is no cost associated with demand shifting.
In this context, the same PAR reduction $(30\%)$ can be obtained for all price increments. 
This follows since the lack of shifting cost maximises the amount of shifted demand.
To illustrate this scenario in more detail, Fig. \ref{fig:simul2} provides the hourly demand of cluster 1 and  the  aggregated  load  profile  experienced  by  the  grid. The observed load profiles can be perceived as unrealistic because the behavior of the consumers is not reflected by the economic signals embedded in the price structure. This result highlights the importance of the incorporation of the shifting cost in \eqref{eq:LLs} to realistically model the demand response of consumers.

\subsection{{Comparison between the IBP and ToU schemes}}
\label{subsec:comparison}
\begin{figure}
    \centering
     \includegraphics[width=3.6in]{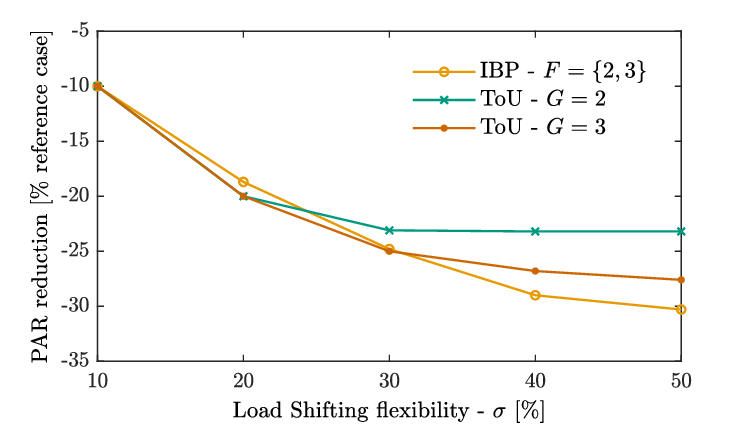}
     \vspace{-5mm}
    \caption{{Optimal PAR reduction versus  load shifting flexibility $\sigma$ for IBP and ToU schemes with two and three consumption-blocks and price tiers respectively.}}
    \label{fig:comparepar}
    \vspace{-3mm}
\end{figure}

{
In this section we compare the performance of our proposed IBP scheme with a state-of-the-art scheme adopted in the literature. 
In particular, we compare the results of our proposed solution against a ToU pricing scheme, as both schemes can potentially promote demand response among residential consumers while providing a long-term price structure. To enable an unbiased comparison,  we designed an optimal price structure for the ToU scheme by forming a suitable bilevel problem, similar to the one considered for the IBP scheme.
The analytical approach for optimizing ToU prices is presented in Appendix B.

We first evaluate the PAR reduction for different IBP and ToU price structures for five scenarios of load shifting flexibility $\sigma=\{10, 20, 30, 40, 50\}$\%.
For the IBP scheme, we considered the cases of two and three consumption-blocks ($F=\{2,3\}$) and for the ToU scheme we considered   two and three price tiers ($G=\{2,3\}$). In all cases, the implemented shifting cost was $\tau=0.03$ \pounds/kWh$^2$.
The resulting PAR reduction, as a percentage of the reference case, is depicted in Fig. \ref{fig:comparepar}.
From Fig.  \ref{fig:comparepar}, it follows that the four cases yield the same PAR reduction when the shifting flexibility is low ($\sigma=10\%$).
For $\sigma=20\%$, the ToU schemes perform slightly better (approximately 6\%) than the IBP schemes.
However, as the shifting flexibility increases, the IBP schemes enable larger PAR reduction than the ToU schemes.
For instance, when $\sigma=50\%$ the IBP schemes perform $30\%$ and $10\%$ better than the ToU schemes with $G=2$ and $G=3$ respectively. 
The latter demonstrates that the proposed IBP scheme can effectively accommodate higher load shifting flexibility. 

}

\begin{figure}
    \centering
     \includegraphics[width=3.3in]{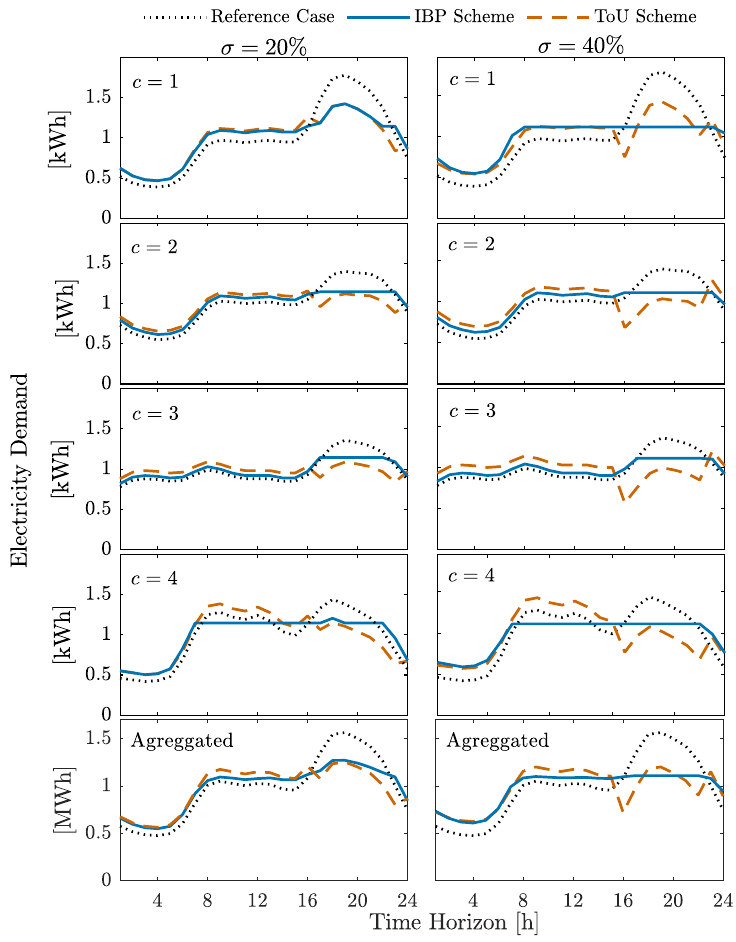}
     \vspace{-3mm}
    \caption{{Electricity demand profiles when load shifting  flexibility is $\sigma=20\%$ (left column) and $\sigma=40\%$ (right column). The black dotted, blue continuous and red dashed lines depict the reference case, and the IBP  ($F=2$) and ToU ($G=2$) scheme responses respectively.}}
    \label{fig:IBPvsToU}
    \vspace{-3mm}
\end{figure}
{
To illustrate the consumer response due to the IBP and ToU schemes, Fig. \ref{fig:IBPvsToU} provides the hourly demand of each cluster $c\in\mathcal{C}=\{1,2,3,4\}$ and the aggregated load profile experienced by the grid, under two load shifting flexibility scenarios, $\sigma=20\%$ (left column) and $\sigma=40\%$ (right column). 
When $\sigma=20\%$, the optimal price structure for the IBP scheme with two consumption-blocks ($F=2$) is given by the prices $\lambda^s_1=0.08$ \pounds/kWh and  $\lambda^s_2=0.1$ \pounds/kWh and the consumption breakpoint $q_1=1.15$ kWh. 
The optimal ToU scheme with two price tiers ($G=2$) is described by the peak price $\lambda^*_2=0.1$ \pounds/kWh within hours $t\in\{17,\cdots,23\}$ and the off-peak price $\lambda^*_1=0.07$ \pounds/kWh in the remaining hours.
When $\sigma=40\%$, the optimal price structure for the IBP scheme with two consumption-blocks ($F=2$) is given by the prices $\lambda^s_1=0.08$ \pounds/kWh and  $\lambda^s_2=0.11$ \pounds/kWh and the consumption breakpoint  $q_1=1.11$ kWh. 
In analogy, the optimal ToU scheme with two price tiers ($G=2$) is described by the peak price $\lambda^*_2=0.09$ \pounds/kWh within hours $t\in\{16,\cdots,23\}$ and the off-peak price $\lambda^*_1=0.07$ \pounds/kWh in the remaining hours.

From the left column in Fig. \ref{fig:IBPvsToU}, where $\sigma=20\%$, it can be observed that both the IBP and ToU schemes result in a similar PAR reduction for clusters $c\in\{1,2,3\}$ and aggregated demand response. However, for cluster $c=4$, the response from the ToU scheme exhibits a peak  shift from 20:00 to 8:00 o'clock. By contrast, the response from the IBP scheme is nearly flat within the same period. 
These trends are more evident when there is larger load shifting flexibility, which is the case in the right column of Fig. \ref{fig:IBPvsToU}, where $\sigma=40\%$. 
From Fig. \ref{fig:IBPvsToU} it can be observed that the IBP scheme has the potential to flatten the demand response of consumers and yield a smooth aggregated load profile. 
By contrast, when the ToU scheme is implemented,  we observe that consumers tend to shift their peak demand from peak to off-peak price times, resulting to an uneven aggregate demand response.
Hence, the proposed IBP scheme avoids concentration of load shifting which justifies its improved  performance, compared to the ToU scheme, under high load shifting flexibility.
}

\subsection{{Performance of the IBP price structure under parametric uncertainty on the demand response model
}}
\label{subsec:robust}

 {
This section aims to demonstrate the robustness of the IBP scheme under parametric uncertainty on consumers shifting behavior.
In particular, we first consider an estimated set of demand response model parameter values ($\tau'$,$\sigma'$) to design an optimal IBP price structure, with set of parameters, associated with $\Xi_U$, denoted by $\Xi^{IBP}_{(\tau',\sigma')}$.
To evaluate the performance of this design under parametric uncertainty, we consider two  sets of potential true parameter values, $\hat{\tau}^*=\{\tau^*_1,\tau^*_2,\cdots\}$ and $\hat{\sigma}^*=\{\sigma^*_1,\sigma^*_2,\cdots\}$.
We then implement the pricing scheme $\Xi^{IBP}_{(\tau',\sigma')}$ assuming some true values ($\tau^*,\sigma^*) \in  \hat{\tau}^* \times \hat{\sigma}^*$ and compute the resulting PAR, which we denote by PAR$_{(\tau^*,\sigma^*),(\tau',\sigma')}$. In addition, we evaluate the feasibility,  with respect to revenue adequacy and bill protection constraints, of the resulting case, allowing a small $2\%$ margin to include a broader set of solutions. 
Finally, we use the selected values ($\tau^*,\sigma^*$) to design an optimal price structure $\Xi^{IBP}_{(\tau^*,\sigma^*)}$ and denote the resulting PAR by
PAR$_{(\tau^*,\sigma^*)}$. 

The above procedure is repeated for all possible pairs of  ($\tau^*,\sigma^*) \in  \hat{\tau}^* \times \hat{\sigma}^*$. 
To estimate the performance of the IBP price structure $\Xi^{IBP}_{(\tau',\sigma')}$ at each different case of true parameter values, the percentage error between the designed and optimal PARs is calculated as follows, 

\begin{equation}
        \text{error}_{(\tau^*,\sigma^*),(\tau'\sigma')}=\frac{\text{PAR}_{(\tau^*,\sigma^*),(\tau',\sigma')} - \text{PAR}_{(\tau^*,\sigma^*)}}{\text{PAR}_{(\tau^*,\sigma^*)}},
		  \label{eq:new3}
\end{equation}
noting that in all feasible scenarios it follows that PAR$_{(\tau^*,\sigma^*)} \leq~$PAR$_{(\tau^*,\sigma^*),(\tau',\sigma')}$.}

\begin{figure}
    \centering
     \includegraphics[width=3.3in]{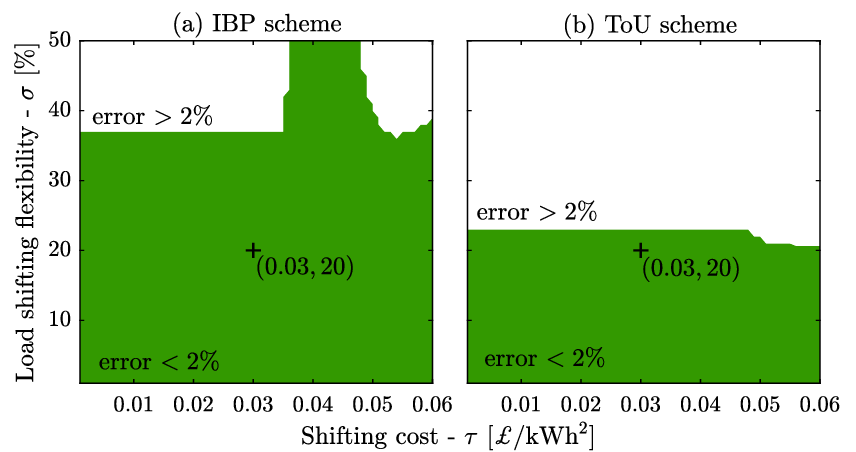}
     \vspace{-3mm}
    \caption{{
    Robustness of the IBP (left) and ToU (right) price schemes to parametric uncertainty. 
    The green regions present the ranges of values where the optimal IBP and ToU price structures that are designed based on estimated parameters values of ($\tau, \sigma) = (0.03  \text{ \pounds/kWh}^2,20\%)$ (pointed with dot), yield a PAR that is no more than $2\%$ larger than the PAR from an optimal price structure design based on these values.
    }}
    \label{fig:robust}
    \vspace{-3mm}
\end{figure}

{For our case study, we consider estimated parameter values given by ($\tau', \sigma') = (0.03  \text{ \pounds/kWh}^2,20\%)$ and true parameter scenario sets given by  $\hat{\sigma}^*=\{0,1,2,\cdots,50\}\%$ and $\hat{\tau}^*=\{0.001,0.002,\cdots,0.06\}$ \pounds/kWh$^2$.
We also set the price increment to $\xi=0.03$ \pounds/kWh. 
The above described approach was then applied, and for each set of parameter values, we used \eqref{eq:new3} to calculate the error between the PARs resulting from designing the IBP scheme based on the estimated and the true values. 
The results of this approach are presented in Fig. \ref{fig:robust} (left), which presents in green color the range of values where the error was below $2\%$ and in white the values where the error was greater than $2\%$. In addition, the  mark $(0.03,20)$ indicates the selected estimated scenario.
From Fig. \ref{fig:robust} (left), it follows that the IBP scheme exhibits significant robustness to parameter variations, enabling a close to optimal response when $\tau^* \in [0, 0.06]\text{\pounds/kWh}^2$ ($[-100\%, + 100\%]$ difference from $\tau' = 0.03  \text{ \pounds/kWh}^2$) and $\sigma^* \in [0, 36]\%$   ( $[-100\%, + 80\%]$ difference from $\sigma = 20\%$).
These results demonstrate that the IBP scheme enables  robust performance against parametric uncertainty.
}

\begin{figure}
    \centering
     \includegraphics[width=3.6in]{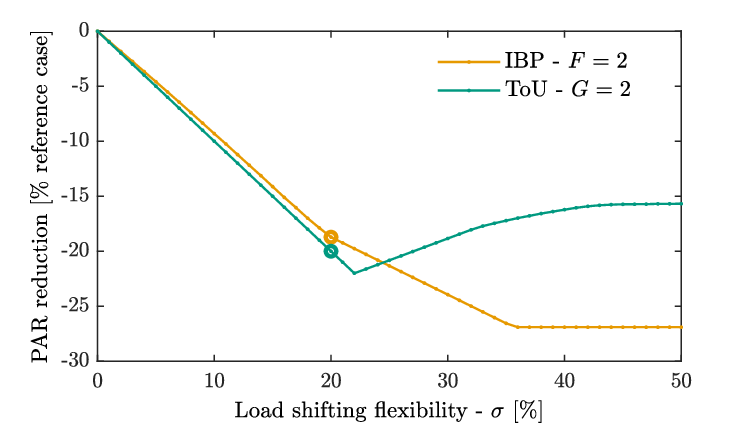}
     \vspace{-3mm}
    \caption{{Performance of the IBP   and ToU schemes under uncertainty in load shifting flexibility.  
    PAR reduction when each scheme is designed based on an estimated value of $\sigma = 20\%$, for a range of true load shifting flexibility values between $0$ and $50\%$.
    }}
    \label{fig:comparrobust}
    \vspace{-3mm}
\end{figure}

{
\subsection{Comparison of the robustness to parametric uncertainty between the IBP and ToU schemes}\label{sec:comparison_rob}

In this section, we repeat  the approach presented in   Section \ref{subsec:robust} on ToU schemes to compare the robustness to parametric uncertainty between the IBP and ToU schemes.
To enable a fair comparison, the same sets of estimated parameter values  and sets of possible true values were considered, i.e. ($\tau', \sigma') = (0.03  \text{ \pounds/kWh}^2,20\%)$,  $\hat{\sigma}^*=\{0,1,2,\cdots,50\}\%$ and $\hat{\tau}^*=\{0.001,0.002,\cdots,0.06\}$ \pounds/kWh$^2$. 
The results of this approach are presented in Fig. \ref{fig:robust} (right).
Figure \ref{fig:robust} (right) presents in green color the range of values where the error, i.e. the difference between the PARs obtained from optimal designs based on the estimated and true values, was below $2\%$ and in white the values where the error was greater than $2\%$. 
From Fig. \ref{fig:robust}, it follows that the proposed IBP scheme exhibits improved robustness to parametric uncertainty, allowing near optimal response for a broader range of parameter values.

To further compare the performance of the IBP and ToU schemes at the presence of parametric uncertainty, we consider its effect on the PAR reduction achieved from these schemes.
In particular, we consider   the same set of estimated parameter values, i.e. ($\tau', \sigma') = (0.03  \text{ \pounds/kWh}^2,20\%)$, and calculate the PAR reduction (compared to flat pricing) when the true values of load shifting flexibility lie in the set $\hat{\sigma}^*=\{0,1,2,\cdots,50\}\%$ and the true value of $\tau$ is fixed at $0.03  \text{ \pounds/kWh}^2$.
We then repeat the approach presented in Section \ref{subsec:robust} but instead of calculating the error, we only consider the PAR reduction.
The results of this approach are presented in Fig. \ref{fig:comparrobust},
which includes the resulting PAR reduction for each case.
In addition, the circle marks indicate the PARs when true load shifting flexibility is equal to the estimated, and hence optimal.
Although the ToU scheme results in a slightly better PAR reduction when $\sigma = 20\%$, we observe that as the true values of load shifting flexibility increase, the IBP scheme yields a significantly larger PAR reduction compared to the ToU scheme.
More precisely, the ToU scheme yields increasing PAR when $\sigma > 22\%$ while the IBP scheme yields a decreasing PAR for values of $\sigma$ up to $35\%$ and a constant PAR when $\sigma$ increases further.
}

\subsection{{Scalability of the solution algorithm}} 
\label{subsec:scalability}
{
This subsection investigates the performance of the solution algorithm as the number of clusters increases to demonstrate its scalability.
A total of seven cases were considered, with the number of  clusters taking values in $C=\{3,5,\cdots,15\}$.  
In addition, the hourly baseline demand $D_{tc}$ and the number of consumers per cluster $n_c$ for each $c\in\mathcal{C}$ were  randomly set. 
Each scenario was simulated 20 times with simulation parameters set at two consumption-blocks ($F=2$), pricing increment $\xi=0.03$ \pounds/kWh, penetration of flexible loads $\sigma=5\%$ and shifting cost $\tau=0.03$ \pounds/kWh$^2$.
}
\begin{figure}
    \centering
     \includegraphics[width=3.5in]{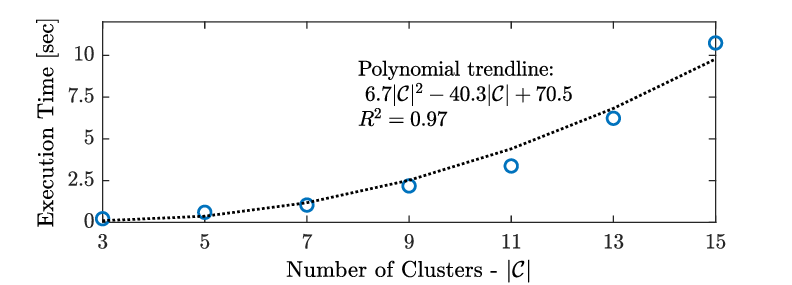}
     \vspace{-5mm}
    \caption{{Time required to solve the MILP problem versus number of clusters.}}
    \label{fig:scalability}
    \vspace{-3mm}
\end{figure}
{
Figure \ref{fig:scalability} provides the average run time value for each scenario. 
From Figure \ref{fig:scalability}, it follows that the time required for the solution algorithm to obtain a solution is proportional to the square of the number of clusters. 
The latter demonstrates the scalability of the proposed approach in regards to the number of clusters.
}

\section{Conclusions}
\label{sec:Conclusions}

We considered the problem of minimizing the PAR of the power system with suitable price incentive schemes for consumers.
In accordance to this, the \textit{Intraday Block Pricing} scheme has been proposed.
Under the IBP scheme, electricity charging is based on intraday pricing periods of equal duration, where the consumption in each pricing period spans a number of pre-defined blocks of increasing price.
A bilevel optimization problem has been formulated for the parameter design of the IBP scheme, aiming to minimize the PAR and simultaneously protect the utility revenue and the consumer bills.
A solution approach has been developed to solve the considered bilevel optimization problem, by converting it into an equivalent single-level Mathematical Program  with  Equilibrium  Constraints, which is then relaxed into a Mixed Integer Linear Program. 
In addition, a lower bound to its cost was obtained to evaluate the conservativeness of the solution approach.
The applicability and implementability of the IBP scheme are demonstrated with numerical simulations on a realistic setting, that demonstrate a significant reduction in the PAR, compared to flat pricing, as well as simultaneous benefits for the utility company and the consumer.
Furthermore, the simulation results yielded equivalent solutions for the original and relaxed problems, demonstrating the effectiveness of the proposed solution approach.
{Finally, simulation results demonstrate the robustness of the IBP scheme to parametric uncertainty, and also improved performance compared to an optimized state-of-the-art  scheme (ToU scheme), when load shifting flexibility is high.}

Future research aims to explore optimal pricing strategies in the context of social welfare maximization and the implications on the pricing strategies and profit of the retailer. 
In addition, we aim to investigate the optimal IBP price structure design problem to accommodate further penetration of intermittent renewable generation and plug-in hybrid electric vehicles. 
{Finally, we aim to develop stochastic optimization methodologies to incorporate  uncertainty in the model parameters and wholesale prices  in the design of the IBP prices in an effective and  computationally efficient manner.
}
{

{
\section*{Appendix A: Incorporation of Network Constraints}

This appendix presents how network constraints may be incorporated in the proposed problem formulation, associated with the design of the IBP scheme prices. 
In addition, it discusses how the solution approach may be extended to include this modification.

To define suitable network constraints, we consider a lossless distribution network with radial structure.
In addition, we assume that the distribution network is three-phased balanced and hence  can be represented through an equivalent single-phase model, similar to \cite{yang2018model}.
To define the network constraints, we consider a distribution network with a set of  buses   denoted by $\mathcal{N}$.
In addition,  we denote the maximum power transfer towards each bus $j \in \mathcal{N}$ by $p^{\max}_j$.
Then, the following constraint can be incorporated into the upper-level problem~\eqref{eq:UL}
\begin{equation}\label{constr_line_transfer}
    \sum_{c\in\mathcal{C}}\sum_{f\in\mathcal{F}} n_{c,j} d^s_{tcf} \leq p^{\max}_j, \forall t\in\mathcal{T}, j \in \mathcal{N},
\end{equation}
where $n_{c,j}$ corresponds to the number of consumers in cluster $c$ and bus $j$.

Since \eqref{constr_line_transfer} is a linear constraint, the proposed solution approach, described in Section \ref{sec:solutionApproach}, can be directly extended to include this case.
In particular, this will  incorporate  $|\mathcal{N}| \times |\mathcal{T}|$ linear constraints in the resulting MILP problem, but will not otherwise complicate the solution approach.
In practical terms, it should be noted that such network constraints will  be satisfied under the application of the proposed IBP scheme, provided that they are satisfied under flat pricing.
The latter follows since the  IBP scheme results in a  non-increase of the individual consumer peak demand and the fact that the considered network has a radial structure.
}

{\section*{Appendix B:  ToU Pricing Scheme Optimal Design}}

This appendix presents our approach to obtain an optimal design of prices associated with the ToU scheme.

The ToU price structure includes a set of prices $\lambda^{ToU}_t$ associated with time-slots $t\in\mathcal{T}$.
Following the work presented in \cite{ferreira2013}, we define a limited amount of price tiers $g\in\mathcal{G}$ within the day,
where $\mathcal{G}=\{1,\cdots,G\}$. Each price tier consists of a set of hours for which prices have the same value $\lambda^*_g$. To enhance intuition among consumers, it is also common practice to design a unimodal ToU price profile $\lambda^{ToU}_t, \forall(t\in\mathcal{T})$, i.e. a price profile with a single peak value.

To optimally design the ToU price structure we define a similar problem to Problem \ref{Problem_statement}.
In particular, we adapt the bilevel problem \eqref{eq:UL}, \eqref{eq:LLs} to incorporate the features of the ToU scheme. 
The adapted bilevel problem is presented in \eqref{eq:ToU} below, where the upper level  designs the price structure and the lower level  estimates the demand response of consumers. Note that we use $d^{ToU}_{tc}$, $\forall(t\in\mathcal{T},c\in\mathcal{C})$ to denote the demand under the ToU scheme for cluster $c$ at time $t$. 

\begin{subequations}\allowdisplaybreaks \label{eq:ToU}
    \begin{alignat}{3}
         \min_{\Xi^{ToU}_U}&                 \frac{\displaystyle\max_{t\in\mathcal{T}}\Bigg\{\sum_{c\in\mathcal{C}}n_cd^{ToU}_{tc}\Bigg\}}{\displaystyle\frac{1}{T}\sum_{t\in\mathcal{T}}\sum_{c\in\mathcal{C}}n_cd^{ToU}_{tc}} \span\span \label{eq:tou1}\\
         \text{s.t. }
         & \lambda^{ToU}_t=\sum_{g\in\mathcal{G}}\lambda^*_g\theta_{tg}, 
            &\forall(t\in\mathcal{T}),  \label{eq:tou2}\\
         &\sum_{g\in\mathcal{G}} \theta_{tg}=1,           
            &\forall(t\in\mathcal{T}),  \label{eq:tou3}\\
         & \lambda^{ToU}_{(t+1)}+ \Theta\Omega^+_t\geq\lambda^{ToU}_t, 
            & \forall(t\in\mathcal{T}^*),  \label{eq:tou4}\\
          & \lambda^{ToU}_{(t+1)}+\Theta\Omega^-_t\leq\lambda^{ToU}_t,
            & \forall(t\in\mathcal{T}^*), \label{eq:tou5}\\
         & \sum_{t\in\mathcal{T}} \Omega^+_{t}\leq G-1,\ \sum_{t\in\mathcal{T}} \Omega^-_{t}\leq G-1, \span  \label{eq:tou6}\\
         & \sum_{t\in\mathcal{T}}\sum_{c\in\mathcal{C}}\lambda^{ToU}_tn_cd^{ToU}_{tc}\geq r\sum_{t\in\mathcal{T}}\sum_{c\in\mathcal{C}}\lambda^w_tn_cd^{ToU}_{tc}, \span   \label{eq:tou7}\\
         & \sum_{t\in\mathcal{T}}\lambda^{ToU}_tD_{tc} \leq\lambda^o\sum_{t\in\mathcal{T}}D_{tc},
            & \forall(c\in\mathcal{C}), \label{eq:tou8}\\
         \min_{\Xi^{ToU}_L}&
		  \sum_{c\in\mathcal{C}} n_c \sum_{t\in\mathcal{T}} (\lambda^{ToU}_td^{ToU}_{tc} + \frac{\tau_{c}}{2}{d^{sh}_{tc}}^2), \span \label{eq:tou9}\\
         \text{s.t. }
         & d^{ToU}_{tc}-d^{sh}_{tc}=D_{tc} : \rho_{tc}, 
            & \forall(t\in\mathcal{T},c\in\mathcal{C}),  \label{eq:tou10}\\
        & 0\leq d^{ToU}_{tc}, 
            & \forall(t\in\mathcal{T},c\in\mathcal{C}),  \label{eq:tou11}\\
        & -\sigma_cD_{tc}\leq d^{sh}_{tc}, & \forall(t\in\mathcal{T},c\in\mathcal{C}),  \label{eq:tou12}\\
        & d^{sh}_{tc}\leq \sigma_cD_{tc},  
            & \forall(t\in\mathcal{T},c\in\mathcal{C}),  \label{eq:tou13}\\
        & \sum_{t\in\mathcal{T}}d^{sh}_{tc}=0, 
            & \forall(c\in\mathcal{C})  \label{eq:tou14},
    \end{alignat}
\end{subequations} 
where $\Theta$ is a sufficiently large positive constant.
The set of decision variables for the upper and lower level problems are given by $\Xi^{ToU}_U=\{\lambda^{ToU}_t,\ (\Omega^-_t,\Omega^+_t)\in\{0,1\},\ \forall(t\in\mathcal{T});\ \lambda^*_g,\ \forall(g\in\mathcal{G});\ \theta_{tg}\in\{0,1\},\ \forall(t\in\mathcal{T},g\in\mathcal{G})\}$ and $\Xi^{ToU}_L=\{d^{ToU}_{tc},\ d^{sh}_{tc},\  \forall(t\in\mathcal{T},c\in\mathcal{C})\}\}$, respectively.
We also define the set $\mathcal{T}^*=\{1,\cdots,T-1\}$, where $T=|\mathcal{T}|$. 

The objective function of the upper level problem \eqref{eq:tou1} minimizes the PAR of the power system. 
Constraints \eqref{eq:tou2} and \eqref{eq:tou3} limit the number of price tiers to a predefined value $|\mathcal{G}|$. 
Constraints \eqref{eq:tou4}-\eqref{eq:tou6} ensure that the resulting price profile is unimodal. 
Revenue adequacy and consumer bill protection are ensured in \eqref{eq:tou7} and \eqref{eq:tou8}, respectively. The objective function of the lower level problem \eqref{eq:tou9} minimizes the cost of consumers and constraint \eqref{eq:tou10} ensures demand balance. Finally constraints \eqref{eq:tou12}-\eqref{eq:tou14} describe the load shifting flexibility of consumers.
The solution approach for problem \eqref{eq:ToU} is similar to that presented in Section \ref{sec:solutionApproach}, and hence omitted.
}

\bibliographystyle{IEEEtran}      
\bibliography{ref}
\let\mybibitem\bibitem

\end{document}